\pdfoutput=1

\documentclass{article}

\usepackage{url} %for formatting urls
\usepackage{array} %for tabular fixed-width columns
\usepackage{amsthm, amsmath, graphicx, amsfonts, mathrsfs, float, enumerate, multirow, rotating}
\usepackage{longtable, subcaption}
\usepackage{array}
\usepackage{multicol,multirow}
\usepackage[labelfont=bf]{caption}
\usepackage[utf8]{inputenc}
\usepackage[english]{babel}
\usepackage[style=alphabetic]{biblatex}
\usepackage{subcaption}
\DeclareCaptionSubType*[Alph]{table}
\DeclareCaptionLabelFormat{mystyle}{Table~\bothIfFirst{#1}{ }#2}
\usepackage[utf8]{inputenc}
\usepackage{authblk}
\usepackage[margin=1in]{geometry}

\newcolumntype{M}[1]{>{\centering\arraybackslash}m{#1}}
\newcolumntype{N}{@{}m{0pt}@{}}

\usepackage{lineno,hyperref}

\modulolinenumbers[1]

\usepackage{graphicx}

\usepackage{amssymb}

\usepackage{amsthm}

\usepackage{lineno}

\usepackage{hyperref}
\usepackage{xcolor}
\pdfoutput=1

\usepackage{multirow}

\usepackage{amsmath,amsfonts,amssymb,amsthm,epsfig,epstopdf,url,array}
\theoremstyle{plain}
\newtheorem{theorem}{Theorem}

\newtheorem{lemma}[theorem]{Lemma}

\theoremstyle{definition}
\newtheorem{definition}{Definition}

\allowdisplaybreaks
\newtheorem*{remark}{Remark}

\begin{document}
\title{New approximate-analytical solutions for the nonlinear fractional Schr\"{o}dinger equation with second-order spatio-temporal dispersion via double Laplace transform method}

\author[1]{Mohammed K. A. Kaabar\thanks{Corresponding author e-mail:mohammed.kaabar@wsu.edu}}
\author[2]{Francisco Mart\'{i}nez}
\author[3]{Jos\'{e} Francisco G\'{o}mez-Aguilar}
\author[4]{Behzad Ghanbari}
\author[5]{Melike Kaplan}
\affil[1]{Department of Mathematics and Statistics, Washington State University, Pullman, WA, USA}
\affil[2]{Department of Applied Mathematics and Statistics, Technological University of Cartagena, Cartagena, Spain}
\affil[3]{CONACyT-Tecnol\'{o}gico Nacional de M\'{e}xico/CENIDET, Interior Internado Palmira S/N, Col. Palmira, C.P. 62490, Cuernavaca Morelos, Mexico}
\affil[4]{Department of Engineering Science, Kermanshah University of Technology, Kermanshah, Iran}
\affil[5]{Department of Mathematics, Art-Science Faculty, Kastamonu University, Kastamonu, Turkey}
\renewcommand\Authands{ and }

\maketitle

\begin{abstract}
  In this paper, a modified nonlinear Schr\"{o}dinger equation with spatio-temporal dispersion is formulated in the senses of Caputo fractional derivative and conformable derivative. A new generalized double Laplace transform coupled with Adomian decomposition method has been defined and applied to solve the newly formulated nonlinear Schr\"{o}dinger equation with spatio-temporal dispersion. The approximate analytical solutions using the proposed generalized method in the sense of Caputo fractional derivative and conformable derivatives are obtained and compared with each other graphically.
\end{abstract}

\section{Introduction}
Fractional derivatives (FDs) which are generalized forms of fractional calculus (FC) have been recently applied in various modeling scenarios arising from phenomena in science and engineering. FDs are used particularly in modeling various complex engineering systems in mechanics. Various systems in science and engineering have been investigated in \cite{ATT, AG, YZ, AMDT, MK, SSR} using analytical and approximate analytical methods to find solutions to the fractional equations in the sense of fractional derivatives and fractional integrals such as Caputo, Riemann-Liouville, and Grünwald–Letnikov \cite{Ekm}. All fractional derivatives such as Caputo, Riemann-Liouville, and Grünwald–Letnikov have shown the property of linearity, but the usual derivative properties such as the derivatives of constant, quotient rule, product rule, and chain rule can not be satisfied using any of those fractional derivatives \cite{AKH}.

Khalil et al. introduced a new generalized local fractional derivative known as conformable derivative which is basically an extension of the usual limit-based derivative \cite{KHH, MLT}. Conformable derivative (CD) occurs naturally and satisfies all properties of usual derivatives \cite{KHH, MLT}. The first author discussed in \cite{mkaabar} three different analytical methods for solving two-dimensional wave equation involving conformable derivative. Introducing CD formulations to (linear/nonlinear) partial differential equations can easily transform them into simpler versions \cite{MLT} than other commonly used fractional derivatives formulations such as Riemann-Liouville, Caputo, and Grünwald–Letnikov because the difficulty of finding analytical solutions using those fractional definitions makes CD formulation a good option for certain cases. The physical interpretation for conformable derivative can be simply presented as modified version of the usual derivative in magnitude and direction \cite{SiLa} (see also \cite{coZ}). Khalil et al. described in \cite{gKH} the geometrical meaning of conformable derivative using the concept of fractional cords. For a comparative review and analysis of all definitions of fractional derivatives and fractional operators, we refer to \cite{mRev23}.

From this newly defined derivative, many essential elements of the mathematical analysis of functions of a real variable have been successfully developed, among which we can mention: Mean Value Theorem, Rolle's Theorem, chain rule, conformable integration by parts formulas, conformable power series expansion and conformable single and double Laplace transform definitions, \cite{KHH, AJ, TO, mkaabar}. The conformable partial derivative of the order $\gamma \in (0,1]$ of the real-valued functions of several variables and conformable gradient vector are defined, and a conformable Clairaut's Theorem for partial derivatives in the conformable sense is proven in \cite{ATN}. In \cite{YzzY}, the conformable Jacobian matrix is introduced; chain rule for multivariable conformable derivative is defined; the relation between conformable Jacobian matrix and conformable partial derivatives is investigated. In \cite{MarT}, two new results on homogeneous functions involving their conformable partial derivatives are introduced, specifically, the homogeneity of the conformable partial derivatives of a homogeneous function and conformable Euler's Theorem.

The partial differential equations (PDEs), particularly the Schr\"{o}dinger equation, have been applied in several applications in physics and engineering due to the importance of this equation in nonlinear optics which can successfully explain the dynamics of optical soliton propagation in optical fibers. Various fractional formulations have been introduced in \cite{YmMwRm, JFmWO} to obtain exact optical soliton solutions for modified nonlinear Schr\"{o}dinger equation (MNLSE) with spatio-temporal dispersion. Finding analytical and approximate analytical solutions for the modified forms of nonlinear fractional Schr\"{o}dinger equation have become a common research interest for physicists and applied mathematicians in the field of optical soliton propagation because of the applications of this equation in plasma, optics, electromagnetism, fluid dynamics, and optical communication \cite{JFmWO, YmMwRm}. The dynamics of optical soliton propagation in optical fiber can be interpreted from the MNLSE with second-order spatio-temporal dispersion and group velocity dispersion coefficients \cite{JFmWO}. Given $\Psi(x,t)$ as a complex-valued wave function that represents the macroscopic property of wave profile of the spatial and temporal variables which are expressed as $x$ and $t$, respectively. Then, MNLSE can be written as \cite{JFmWO}:\begin{equation} \begin{split}
& i\displaystyle \left(\frac{\partial \Psi}{\partial x}+\omega_{1}\frac{\partial \Psi}{\partial t}\right)+\omega_{2}\frac{\partial^2 \Psi}{\partial t^2}+\omega_{3}\frac{\partial^2 \Psi}{\partial x^2}+\displaystyle \left|\Psi \right|^2 \Psi=0, \\
       \text{where \cite{JFmWO} } & \text{$\omega_{1}$ is proportional to the ratio of group speed};\\
                                  & \text{$\omega_{2}$ is a group velocity dispersion coefficient};\\
                                  & \text{$\omega_{3}$ is a spatial dispersion coefficient};
\end{split}
\end{equation}

To formulate (1) in the sense of fractional derivatives, let us first define Caputo fractional derivative as follows:
\begin{definition}
For $\xi, \gamma > 0$, given two functions: $h(x)$ and $h(t)$ such that for $x,t>0$, the Caputo fractional derivative (CpFD) of $h$ of order $\xi$ and $\gamma$, denoted by $D_{x}^{\xi}(h)(x)$ and $D_{t}^{\gamma}(h)(t)$, respectively where $D_{x}^{\xi}$ and $D_{t}^{\gamma}$ are Caputo derivative operators which can be simply expressed as \cite{CHinD}:\begin{equation}
 \displaystyle D^{\xi}_{x}h(x) = \frac{1}{\xi(\Omega-\xi)}\int_{0}^{x} (x-\eta)^{\Omega-\xi-1}h^{(\Omega)}(\eta)\mathrm{d}\eta;\, \, \text{$\Omega-1 < \xi \leq \Omega$ for $\Omega\in \mathbb{N}$,}
\end{equation}
\begin{equation}
 \displaystyle D^{\gamma}_{t}h(t) = \frac{1}{\Gamma(w-\gamma)}\int_{0}^{t} (t-\mu)^{w-\gamma-1}h^{(w)}(\mu)\mathrm{d}\mu;\, \, \text{$w-1 < \gamma \leq w$ for $w\in \mathbb{N}$.}
\end{equation}
\end{definition}
If $\xi=\Omega$ and $\gamma=w$ where $\Omega, w\in \mathbb{N}$, then $D^{\xi}_{x}h(x) = \frac{d^{\Omega}}{dx^{\Omega}}h(x)$ and $D^{\gamma}_{t}h(t) = \frac{d^{w}}{dt^{w}}h(t)$. CpFD is very useful in science  and engineering due to their important properties such as the inclusion of initial and boundary conditions in the fractional formulation of CpFD \cite{ZCpt, ATT}. Let us now define the Mittag-Leffler function:\begin{definition}
The Mittag-Leffler function, denoted by $E_{\xi,\zeta}(t)$, can be expressed as follows \cite{RdGw}:
\begin{equation}
 \displaystyle E_{\xi,\zeta}(t)= \sum^{\infty}_{m=0} \frac{t^{m}}{\Gamma(\xi m+\zeta)}, \text{where $t,\zeta \in \mathbb{C}$ and $\Re(\xi)>0$}.
\end{equation}
\end{definition} From the above definition, the Mittag-Leffler function, denoted by $E(t,h,c)$, can be written \cite{CHinD} as: $E(t,h,c)=t^{h}E_{1,h+1}(ct)$, and the fractional derivative of Mittag-Leffler function can also be expressed \cite{CHinD} as: $\frac{\partial^{\delta}}{\partial t^{\delta}}(t^{\zeta-1}E_{\xi,\zeta}(ct^{\xi}))=t^{\zeta-\delta-1}E_{\xi,\zeta-\delta}(ct^{\xi})$ where $\delta\geq0$.
let us now define the conformable derivative as follows:
\begin{definition}
Given a function $h:[0,\infty)\rightarrow\Re$ such that for all $t>0$, the conformable derivative (CD) of order $\gamma\in(0,1]$ of $h$, denoted by $M_{\gamma}(h)(t)$, can be represented as:\begin{equation}
 \displaystyle h^{(\gamma)}(t) = M_{\gamma}(h)(t) = \lim_{\xi\to 0} \frac{h(t+\xi t^{1-\gamma}) - h(t)}{\xi}.
\end{equation}
Suppose that $h$ is $\gamma$-differentiable in some $(0,d)$, $d>0$, and the limit of $h^{(\gamma)}(t)$ exists as $t\longrightarrow 0^{+}$, then from CD definition, the following is obtained:
\begin{equation}
 \displaystyle h^{(\gamma)}(0) = M_{\gamma}(h)(0) = \lim_{t \to 0^{+}}h^{(\gamma)}(t).
\end{equation}
\end{definition}
It is also important to define here the conformable integral (ComI) \cite{SiLa} as follows:
\begin{definition}
For $\gamma\in(0,1]$, given a function $h:[0,\infty)\rightarrow\Re$ such that for all $t\geq0$, the $\gamma$th order ComI of $h$ from $0$ to $t$ can be expressed as:\begin{equation}
\begin{split}
 \displaystyle I_{\gamma}(h)(t) = \int_{0}^{^{t}}h(\psi)d_{\gamma}\psi=\int_{0}^{^{t}}h(\psi)\psi^{\gamma -1}d\psi
 \end{split}
\end{equation}
\end{definition}

If we suppose $\gamma=1$, then we have $I_{\gamma}(h)(t)=I_{\gamma=1}(t^{\gamma -1}h)(t)$ which represents the classical improper Riemann integral of a function $h(t)$. Given a continuous function, $h$, on $(0,\infty)$ and for $\gamma\in(0,1)$, then $M_{\gamma}(h)(t)\displaystyle \left[I_{\gamma}(h)(t)\right]=h(t)$.

\begin{lemma} \cite{AJ}
Given a function $h:(c,d)\rightarrow\Re$ as a differentiable function and $\gamma \in (0,1]$. Then for all $c>0$, we have the following:
\begin{equation*}
 \displaystyle I_{\gamma}^{c}M_{\gamma}^{c}(h)(t)=h(t)-h(c).
\end{equation*}
\end{lemma}

In addition, the following theorem \cite{KHH, AKH} shows that $M_{\gamma}$ satisfies all the standard properties of basic limit-based derivative as follows:

\begin{theorem}
For $\gamma\in(0,1]$, given two functions say: $h$ and $v$ to be assumed $\gamma$-differentiable at a point $t$, then the following is obtained:
\begin{itemize}
\item[(1)]   $M_{\gamma}(ch+ev)=cM_{\gamma}(h)+eM_{\gamma}(v)$, for all $c,e\in\Re$.
\item[(2)]   $M_{\gamma}(t^{r})=rt^{r-\gamma}$, for all $r\in\Re$.
\item[(3)]   $M_{\gamma}(hv)=hM_{\gamma}(v)+vM_{\gamma}(h)$.
\item[(4)]   $M_{\gamma}(\frac{h}{v})=\frac{vM_{\gamma}(h)-hM_{\gamma}(v)}{v^{2}}$.
\item[(5)]   $M_{\gamma}(\lambda)=0$, for all constant functions $h(t)=\lambda$.
\item[(6)]   If $h$ is assumed to be a differentiable function, then $M_{\gamma}(h)(t)=t^{1-\gamma}\frac{dh}{dt}$.
\end{itemize}
\end{theorem}

The conformable partial derivative of a real valued function with several variables is defined in
\cite{ATN, YzzY} as follows:

\begin{definition}
Let $h$ be a real-valued function with $n$ variables and $\textbf{c}=(c_{1},c_{2},...,c_{n})\in \Re^{n}$ be a
point whose $i^{th}$ component is positive. Then, we have:
\begin{equation*}
 \displaystyle \lim_{\xi\to 0} \frac{h(c_{1},c_{2},...,c_{i}+\xi c_{i}^{1-\gamma},...,c_{n})-h(c_{1},c_{2},...,c_{n})}{\xi}.
\end{equation*}
If the limit exists, the $i^{it}$ conformable partial derivative of $h$ of the order $\gamma \in(0,1]$ at
$\textbf{c}$ is denoted by $\frac{\partial^{\gamma}}{\partial x_{i}^{\gamma}}h(\textbf{c})$.
\end{definition}

\begin{remark}
Let $\gamma \in (0,\frac{1}{\Lambda}]$, $\Lambda \in \mathbb{Z^{+}}$ and $h$ be a real-valued function with $n$ variables defined on an open set $D$, such that for all $(x_{1}, x_{2},..., x_{n})\in D$, each $x_{i}>0$. The function, $h$, is said to be $C_{\gamma}^{\Lambda}(D,\Re)$ if all its conformable partial derivatives of order less than or equal to $\Lambda$ exist and are continuous on $D$, \cite{MarT}.
\end{remark}

For more information about other related fractional derivatives' definitions and their physical and geometrical interpretations, we refer to \cite{mkaabar}.The main goal of this paper is to obtain approximate-analytical solutions for MNLSE in (1) using double Laplace transform method in the sense of Caputo and conformable derivatives. From definition 3, Let us first formulate the MNLSE in (1) in the sense of CD as follows:\begin{equation}
\begin{split}
& i\displaystyle \left(D_{x}^{\gamma}\Psi \displaystyle \left(\frac{x^{\gamma}}{\gamma},\frac{t^{\delta}}{\delta}\right)+\omega_{1}D_{t}^{\delta}\Psi\displaystyle \left(\frac{x^{\gamma}}{\gamma},\frac{t^{\delta}}{\delta}\right)\right)+
\omega_{2}D_{t}^{2\delta}\Psi\displaystyle \left(\frac{x^{\gamma}}{\gamma},\frac{t^{\delta}}{\delta}\right)+\omega_{3}D_{x}^{2\gamma}\Psi\displaystyle \left(\frac{x^{\gamma}}{\gamma},\frac{t^{\delta}}{\delta}\right) \\
&+\displaystyle \left|\Psi \right|^2 \Psi=0; \\
& iD_{x}^{\gamma}\Psi\displaystyle \left(\frac{x^{\gamma}}{\gamma},\frac{t^{\delta}}{\delta}\right)+\omega_{1}iD_{t}^{\delta}\Psi\displaystyle \left(\frac{x^{\gamma}}{\gamma},\frac{t^{\delta}}{\delta}\right)+
\omega_{2}D_{t}^{2\delta}\Psi\displaystyle \left(\frac{x^{\gamma}}{\gamma},\frac{t^{\delta}}{\delta}\right)+\omega_{3}D_{x}^{2\gamma}\Psi\displaystyle \left(\frac{x^{\gamma}}{\gamma},\frac{t^{\delta}}{\delta}\right) \\
&+\displaystyle \left|\Psi \right|^2 \Psi=0; \\
& \text{where $i=\sqrt{-1}$, $0< \gamma,\delta \leq 1$, $t,x>0$.} \\
\end{split}
\end{equation}

From definition 1, MNLSE in (1) can be similarly formulated in the sense of CpFD as follows:\begin{equation} \begin{split}
& i\displaystyle \left(\frac{\partial^{\gamma} \Psi(x,t)}{\partial x^{\gamma}}+\omega_{1}\frac{\partial^{\delta} \Psi(x,t)}{\partial t^{\delta}}\right)+\omega_{2}\frac{\partial^{2\delta} \Psi(x,t)}{\partial t^{2\delta}}+\omega_{3}\frac{\partial^{2\gamma} \Psi(x,t)}{\partial x^{2\gamma}}+\displaystyle \left|\Psi \right|^2 \Psi=0; \\
& i\frac{\partial^{\gamma} \Psi(x,t)}{\partial x^{\gamma}}+\omega_{1}i\frac{\partial^{\delta} \Psi(x,t)}{\partial t^{\delta}}+\omega_{2}\frac{\partial^{2\delta} \Psi(x,t)}{\partial t^{2\delta}}+\omega_{3}\frac{\partial^{2\gamma} \Psi(x,t)}{\partial x^{2\gamma}}+\displaystyle \left|\Psi \right|^2 \Psi=0; \\
& \text{where $i=\sqrt{-1}$, $0< \gamma,\delta \leq 1$, $t,x>0$.} \\
\end{split}
\end{equation}

\section{The analytical solutions of nonlinear fractional Schr\"{o}dinger equation}

Many differential equations can be easily solved by applying the method of Laplace transform (DT) for a single-variable function. New generalized forms of the classical Laplace transform methods such as double Laplace transform and multiple Laplace transform have been first introduced in \cite{EHLAPLACE} to solve partial differential equations (PDEs). Recently, double Laplace transform (DLTr) has become an interesting topic of research for many mathematicians and researchers \cite{Adam1,Adam2,Adam3} because not many research studies have been done on this topic \cite{LDDL1} and the need to find an efficient method for solving PDEs. Applying the DLTr in the sense of fractional derivatives has been rarely discussed and it is considered as an open problem \cite {LABALN}. DLTr has been successfully introduced in solving some fractional differential equations (FDEs) such as the fractional heat equation and the fractional telegraph equation via the definition of CpFD \cite{RdGw,LABALN}. According to our knowledge, the generalized DLTr method has never been applied before for solving the MNLSE in (1) in the senses of CpFD and CD. Therefore, the results in this work are new and worthy.

The first author has defined in \cite{mkaabar} the conformable double Laplace transform (CmDLTr) as follows:
\begin{definition}
Given a function, $\Psi\displaystyle \left(\frac{x^{\gamma}}{\gamma},\frac{t^{\delta}}{\delta}\right):[0,\infty)\rightarrow\Re$ such that for all $x,t>0$, the CmDLTr of order $\gamma,\delta\in(0,1]$ of $\Psi\displaystyle \left(\frac{x^{\gamma}}{\gamma},\frac{t^{\delta}}{\delta}\right)$, denoted by $\ell_{\gamma\delta}^{xt}\displaystyle \left[ \Psi\displaystyle \left(\frac{x^{\gamma}}{\gamma},\frac{t^{\delta}}{\delta}\right)\right]$, starting from $0$ can be expressed as follows:
\begin{equation}
 \begin{split}
 \displaystyle \ell_{\gamma\delta}^{xt}\displaystyle \left[ \Psi\displaystyle \left(\frac{x^{\gamma}}{\gamma},\frac{t^{\delta}}{\delta}\right)\right]
&=\ell_{\gamma}^{x}\ell_{\delta}^{t}\displaystyle \left[ \Psi\displaystyle \left(\frac{x^{\gamma}}{\gamma},\frac{t^{\delta}}{\delta}\right)\right]
 =\tilde{\Psi}_{\gamma\delta}^{xt}(s_{1},s_{2})\\
&=\int_{0}^{\infty}\int_{0}^{\infty}e^{-(s_{1}\frac{x^{\gamma}}{\gamma}+s_{2}\frac{t^{\delta}}{\delta})}\Psi\displaystyle \left(\frac{x^{\gamma}}{\gamma},\frac{t^{\delta}}{\delta}\right)x^{\gamma -1}t^{\delta -1}\,dx\,dt.
 \end{split}
 \end{equation}
\end{definition}
where $s_{1},s_{1}\in\mathbb{C}$. If the integral in the above definition exists, then this definition holds true.

To define the double Laplace transform in the sense of Caputo partial fractional derivatives, let's assume that $\tilde{\Psi}^{xt}(s_{1},s_{2})=\int_{0}^{\infty}\int_{0}^{\infty}e^{-(s_{1}x+s_{2}t)}\Psi(x,t)\,dx\,dt.$ From  \cite{RdGw}, theorems 3.1 and 3.3 in \cite{LABALN}, and theorem 2 in \cite{Eurp}, the Caputo double Laplace transform can be defined as follows:
\begin{definition}
Given a function, $\Psi(x,t):[0,\infty)\rightarrow\Re$ such that for all $x,t>0$, the double Laplace transform of the Caputo partial fractional derivatives (CpDLTr) of $\Psi(x,t)$ of orders $\xi$ and $\gamma$ where $\xi \in (j-1,j]$ and $\gamma \in (b-1,b]$ such that $\xi, \gamma >0$ and $j,b \in \mathbb{N}$, denoted by $\frac{\partial ^{\xi}}{\partial x^{\xi}}\Psi(x,t)$ and $\frac{\partial ^{\gamma}}{\partial t^{\gamma}}\Psi(x,t)$, respectively can be expressed as:
\begin{equation}
 \displaystyle \ell_{\xi}^{xt} \displaystyle \left[\frac{\partial^{\xi}}{\partial x^{\xi}}\Psi(x,t) \right]=s_{1}^{\xi}\tilde{\Psi}^{xt}_{\xi}(s_{1},s_{2})-\sum_{i=0}^{j -1}s_{1}^{\xi -1-i}\ell_{t}\displaystyle \left[\frac{\partial^{i}\Psi(0,t)}{\partial x^{i}}\right].
\end{equation}
\begin{equation}
 \displaystyle \ell_{\gamma}^{xt} \displaystyle \left[\frac{\partial^{\gamma}}{\partial t^{\gamma}}\Psi(x,t) \right]=s_{2}^{\gamma}\tilde{\Psi}^{xt}_{\gamma}(s_{1},s_{2})-\sum_{a=0}^{b -1}s_{1}^{\gamma -1-a}\ell_{x}\displaystyle \left[\frac{\partial^{a}\Psi(x,0)}{\partial t^{a}}\right].
\end{equation}
\end{definition}

The double Laplace transform in the sense of conformable partial fractional derivatives can be similarly defined \cite{mkaabar} as follows:
\begin{definition}
Given a function, $\Psi\displaystyle \left(\frac{x^{\gamma}}{\gamma},\frac{t^{\delta}}{\delta}\right):[0,\infty)\rightarrow\Re$ such that for all $x,t>0$, the double Laplace transform of the conformable partial fractional derivatives (CmDLTr) of $\Psi\displaystyle \left(\frac{x^{\gamma}}{\gamma},\frac{t^{\delta}}{\delta}\right)$ of orders $\gamma$ and $\delta$ where $\gamma, \delta \in (0,1]$, denoted by $\frac{\partial ^{\gamma}}{\partial x^{\gamma}}\Psi\displaystyle \left(\frac{x^{\gamma}}{\gamma},\frac{t^{\delta}}{\delta}\right)$ and $\frac{\partial ^{\delta}}{\partial t^{\delta}}\Psi\displaystyle \left(\frac{x^{\gamma}}{\gamma},\frac{t^{\delta}}{\delta}\right)$, respectively can be written as:
\begin{equation}
 \displaystyle \ell_{\gamma}^{xt} \displaystyle \left[\frac{\partial^{\gamma}}{\partial x^{\gamma}}\Psi\displaystyle \left(\frac{x^{\gamma}}{\gamma},\frac{t^{\delta}}{\delta}\right) \right]=s_{1}^{\gamma}\tilde{\Psi}^{xt}_{\gamma}(s_{1},s_{2})-\sum_{i=0}^{j -1}s_{1}^{\gamma -1-i}\ell_{t}\displaystyle \left[\frac{\partial^{i}\Psi(0,t)}{\partial x^{i}}\right].
\end{equation}
\begin{equation}
 \displaystyle \ell_{\delta}^{xt} \displaystyle \left[\frac{\partial^{\delta}}{\partial t^{\delta}}\Psi\displaystyle \left(\frac{x^{\gamma}}{\gamma},\frac{t^{\delta}}{\delta}\right) \right]=s_{2}^{\delta}\tilde{\Psi}^{xt}_{\delta}(s_{1},s_{2})-\sum_{a=0}^{b -1}s_{1}^{\delta -1-a}\ell_{x}\displaystyle \left[\frac{\partial^{a}\Psi(x,0)}{\partial t^{a}}\right].
\end{equation}
\end{definition}

The first author has proved the existence and uniqueness of CmDLTr in \cite{mkaabar}, while the existence and uniqueness of CpDLTr have been discussed in \cite{LABALN}. It is obvious that the formulas of the double Laplace transform in definition 7 and definition 8 are the same when $j,b=1$. Therefore, the general definition of CpDLTr coincides with the general defintion of CmDLTr when $j,b=1$. The properties of CmDLTr and CpDLTr have been discussed in \cite{Ozzz} and \cite{Adam4}, respectively. For $\gamma,\delta \in (0,1]$, let us now define the formula of the inverse fractional double Laplace transform \cite{LDDL1,RdGw} for both conformable and Caputo fractional derivatives, denoted by $(\ell_{\gamma\delta}^{xt})^{-1}[\tilde{\Psi}^{xt}_{\gamma\delta}(s_{1},s_{1})]$, as follows:

\begin{definition}
Given an analytic function: $\tilde{\Psi}^{xt}_{\gamma\delta}(s_{1},s_{2})$, for all $s_{1},s_{2}\in\mathbb{C}$ and for $\gamma,\delta \in (0,1]$ such that $Re\{s_{1}\geq\eta\}$ and $Re\{s_{2}\geq\sigma\}$, where $\eta,\sigma\in\Re$, then, the inverse fractional double Laplace transform (IFDLT) can be expressed \cite{mkaabar} as follows:
\begin{equation}
 \begin{split}
 &\displaystyle (\ell_{\gamma\delta}^{xt})^{-1}[\tilde{\Psi}^{xt}_{\gamma\delta}(s_{1},s_{2})]
 =(\ell_{\gamma}^{x})^{-1}(\ell_{\delta}^{t})^{-1}[\tilde{\Psi}^{xt}_{\gamma}(s_{1},s_{2})] \\
 &=\frac{-1}{4 \pi^{2}}\int_{\varrho-i\infty}^{\varrho+i\infty}\int_{\varsigma-i\infty}^{\varsigma+i\infty}e^{s_{1}x}
 e^{s_{2}t}\tilde{\Psi}^{xt}_{\gamma\delta}(s_{1},s_{2})\,ds_{1}\,ds_{2}
 \end{split}
 \end{equation}
\end{definition}

To solve the MNLSE in (1) in the senses of CpFD and CD (see equations (9) and (8)) by the methods of CpDLTr and CmDLTr. respectively, let us first re-write both Equation (9) and Equation (8) as follows:

\begin{equation}
\begin{split}
&\frac{\partial^{2\gamma} \Psi(x,t)}{\partial x^{2\gamma}}=-\frac{\omega_{2}}{\omega_{3}} \frac{\partial^{2\delta} \Psi(x,t)}{\partial t^{2\delta}}-\frac{i}{\omega_{3}} \frac{\partial^{\gamma} \Psi(x,t)}{\partial x^{\gamma}}-\frac{\omega_{1}}{\omega_{3}} i\frac{\partial^{\delta} \Psi(x,t)}{\partial t^{\delta}}- \frac{1}{\omega_{3}} \displaystyle \left|\Psi \right|^2 \Psi. \\
&\text{subject to the following initial and boundary conditions:} \\
&\text{$\Psi(x,0)=a_{0}(x)$ and $\frac{\partial\Psi(x,0)}{\partial t}=a_{1}(x)$}. \\
&\text{$\Psi(0,t)=b_{0}(t)$ and $\frac{\partial\Psi(0,t)}{\partial x}=b_{1}(t)$}. \\
&\text{where $i=\sqrt{-1}$, $0< \gamma,\delta \leq 1$, $t,x>0$; $x,t \in\Re^{+}$, and $a_{0},a_{1},b_{0},b_{1} \in \mathbb{C}(\Re^{+},\Re^{+})$.}\\
\end{split}
\end{equation}

\begin{equation}
\begin{split}
&D_{x}^{2\gamma}\Psi\displaystyle \left(\frac{x^{\gamma}}{\gamma},\frac{t^{\delta}}{\delta}\right)=-\frac{\omega_{2}}{\omega_{3}} D_{t}^{2\delta}\Psi\displaystyle \left(\frac{x^{\gamma}}{\gamma},\frac{t^{\delta}}{\delta}\right)-\frac{i}{\omega_{3}} D_{x}^{\gamma}\Psi\displaystyle \left(\frac{x^{\gamma}}{\gamma},\frac{t^{\delta}}{\delta}\right)-\frac{\omega_{1}} {\omega_{3}}i D_{t}^{\delta}\Psi\displaystyle \left(\frac{x^{\gamma}}{\gamma},\frac{t^{\delta}}{\delta}\right) \\
&- \frac{1}{\omega_{3}} \displaystyle \left|\Psi \right|^2 \Psi. \\
&\text{subject to the following initial and boundary conditions:} \\
&\text{$\Psi\displaystyle \left(\frac{x^{\gamma}}{\gamma},0\right)=n_{0}\displaystyle \left(\frac{x^{\gamma}}{\gamma}\right)$ and $D_{t}\Psi\displaystyle \left(\frac{x^{\gamma}}{\gamma},0\right)=n_{1}\displaystyle \left(\frac{x^{\gamma}}{\gamma}\right)$.} \\
&\text{$\Psi\displaystyle \left(0,\frac{t^{\delta}}{\delta}\right)=m_{0}\displaystyle \left(\frac{t^{\delta}}{\delta}\right)$ and $D_{x}\Psi\displaystyle \left(0,\frac{t^{\delta}}{\delta}\right)=m_{1}\displaystyle \left(\frac{t^{\delta}}{\delta}\right)$.} \\
&\text{where $i=\sqrt{-1}$, $0< \gamma,\delta \leq 1$, $t,x>0$; $x,t \in\Re^{+}$, and $n_{0},n_{1},m_{0}, m_{1} \in \mathbb{C}(\Re^{+},\Re^{+})$.}
\end{split}
\end{equation}

By applying the single Laplace transform to initial and boundary conditions in (16) and (17), respectively, we obtain the following:
\begin{equation}
\begin{split}
&\ell[\Psi(x,0)]=\ell[a_{0}(x)]=\tilde{a_{0}}(s_{1});\ell\displaystyle \left[\frac{\partial\Psi(x,0)}{\partial t}\right]=\tilde{a_{1}}(s_{1}). \\
&\ell[\Psi(0,t)]=\ell[b_{0}(t)]=\tilde{b_{0}}(s_{2});\ell\displaystyle \left[\frac{\partial\Psi(0,t)}{\partial x}\right]=\tilde{b_{1}}(s_{2}).
\end{split}
\end{equation}

\begin{equation}
\begin{split}
&\ell[\Psi(x,0)]=\ell[n_{0}(x)]=\tilde{n_{0}}(s_{1});\ell\displaystyle \left[D_{t}\Psi\displaystyle \left(\frac{x^{\gamma}}{\gamma},0\right)\right]=\ell \displaystyle \left[n_{1}\displaystyle \left(\frac{x^{\gamma}}{\gamma}\right)\right]=\tilde{n_{1}}(s_{1}). \\
&\ell[\Psi(0,t)]=\ell[m_{0}(t)]=\tilde{m_{0}}(s_{2});\ell\displaystyle \left[D_{x}\Psi\displaystyle \left(0,\frac{t^{\delta}}{\delta}\right)\right]=\ell \displaystyle \left[m_{1}\displaystyle \left(\frac{t^{\delta}}{\delta}\right)\right]=\tilde{m_{1}}(s_{2}).
\end{split}
\end{equation}

Let us now apply the CpDLTr (definition 7) to both left-hand and right-hand sides of Equation (16), we obtain:

\begin{equation}
\begin{split}
&\displaystyle \ell^{x}\ell^{t}\displaystyle \left[\frac{\partial^{2\gamma} \Psi(x,t)}{\partial x^{2\gamma}}\right]=
\tilde{\Psi}(s_{1},s_{2})=\frac{s_{1}^{2 \gamma -1}}{s_{1}^{2 \gamma}}\ell^{x}\ell^{t}[b_{0}(t)]+\frac{s_{1}^{s \gamma -2}}{s_{1}^{2 \gamma}}\ell^{x}\ell^{t}[b_{1}(t)] \\
&-\frac{1}{s_{1}^{2 \gamma}}\displaystyle \left[\ell^{x}\ell^{t} \displaystyle \left[\frac{\omega_{2}}{\omega_{3}} \frac{\partial^{2 \delta} \Psi(x,t)}{\partial t^{2 \delta}}+\frac{i}{\omega_{3}} \frac{\partial^{\gamma}\Psi(x,t)}{\partial x^{\gamma}}+i\frac{\omega_{1}}{\omega_{3}}\frac{\partial^{\delta} \Psi(x,t)}{\partial t^{\delta}}\right]+\ell^{x}\ell^{t}\displaystyle \left[\frac{1}{\omega_{3}}\displaystyle \left|\Psi \right|^2 \Psi\right]\right]. \\
&\text{By simplifying the above, we obtain:} \\
&\displaystyle \ell^{x}\ell^{t}\displaystyle \left[\frac{\partial^{2\gamma} \Psi(x,t)}{\partial x^{2\gamma}}\right]=
\tilde{\Psi}(s_{1},s_{2})=\frac{1}{s_{1}}\ell^{x}\ell^{t}[b_{0}(t)]+\frac{1}{s_{1}^{2}}\ell^{x}\ell^{t}[b_{1}(t)] \\
&-\frac{1}{s_{1}^{2 \gamma}}\displaystyle \left[\ell^{x}\ell^{t} \displaystyle \left[\frac{\omega_{2}}{\omega_{3}} \frac{\partial^{2 \delta} \Psi(x,t)}{\partial t^{2 \delta}}+\frac{i}{\omega_{3}} \frac{\partial^{\gamma}\Psi(x,t)}{\partial x^{\gamma}}+i\frac{\omega_{1}}{\omega_{3}}\frac{\partial^{\delta} \Psi(x,t)}{\partial t^{\delta}}\right]+\ell^{x}\ell^{t}\displaystyle \left[\frac{1}{\omega_{3}}\displaystyle \left[\Psi^2 \Psi^{*}\right]\right]\right].\\
&\text{where $\displaystyle \left|\Psi \right|^2 \Psi=\Psi^2 \Psi^{*}$ such that $\Psi^{*}$ is the conjugate of $\Psi$.}
\end{split}
\end{equation}

From the Adomian decomposition method (ADcM) (see \cite{Kunm, Pal} for more information about ADcM), Equation (20) is written according to the following standard operator form for nonlinear partial differential equations (NPDEs):
$L\Psi(x,t)+R\Psi(x,t)+N\Psi(x,t)=S(x,t)$ where $N$ represents the nonlinear differential operator, $L$ represents the $2$nd-order partial differential operator, $R$ represents the remaining linear operator, and $S(x,t)$ represents a source term. This method was first introduced by G. Adomian in 1980s where ADcM is well-known for obtaining a series solution whose each term is obtained recursively \cite{Pal}. Therefore, by applying the method of CpDLTr coupled with ADcM, the decomposition infinite series can be expressed for both linear and nonlinear terms in Equation (2), respectively, as follows:

\begin{equation}
\begin{split}
&\Psi(x,t)=\sum_{i=0}^{\infty}\Psi_{i}(x,t) \\
&N(\Psi(x,t))=\sum_{i=0}^{\infty}\phi_{i}(\Psi(x,t)).
\end{split}
\end{equation}
\begin{equation*}
\begin{split}
&\text{where the above nonlinear term, denoted by $N(\Psi(x,t))$, is represented by infinite} \\
&\text{series of the Adomian polynomials, donated by $\phi_{i}$, which can be expressed \cite{Pal}} \\
&\text{as follows:}
\end{split}
\end{equation*}
\begin{equation}
\begin{split}
&\phi_{i}=\frac{1}{i!}\frac{d^{i}}{d\Omega^{i}} \displaystyle \left[N \displaystyle \left(\sum_{j=0}^{\infty}\Omega^{j}\phi_{j}\right)\right]_{\Omega=0}, i=0,1,2,3,... \\
&\text{so, we can write some of those terms as follows:} \\
&\text{$\phi_{0}=N(\Psi_{0})$; $\phi_{1}=\Psi_{1}N'(\Psi_{0})$; $\phi_{2}=\Psi_{2}N'(\Psi_{0})+\frac{1}{2!}\Psi_{1}^{2}N''(\Psi_{0})$.}
\end{split}
\end{equation}

By applying the standard NPDEs operator form and (21) to Equation (20), we obtain the following:
\begin{equation}
\begin{split}
&\displaystyle \ell^{x}\ell^{t}\displaystyle \left[\sum_{i=0}^{\infty}\Psi_{i}(x,t)\right]=
\frac{1}{s_{1}}\ell^{x}\ell^{t}[b_{0}(t)]+\frac{1}{s_{1}^{2}}\ell^{x}\ell^{t}[b_{1}(t)] \\
&-\frac{1}{s_{1}^{2 \gamma}}\displaystyle \left[\ell^{x}\ell^{t} \displaystyle \left[ R \displaystyle \left[\sum_{i=0}^{\infty}\Psi_{i}(x,t)\right]\right]+\ell^{x}\ell^{t}\displaystyle \left[\sum_{i=0}^{\infty}\phi_{i}    \right]\right]. \\
&\text{where $R[\Psi(x,t)]=\frac{\omega_{2}}{\omega_{3}} \frac{\partial^{2 \delta} \Psi(x,t)}{\partial t^{2 \delta}}+\frac{i}{\omega_{3}} \frac{\partial^{\gamma}\Psi(x,t)}{\partial x^{\gamma}}+i\frac{\omega_{1}}{\omega_{3}}\frac{\partial^{\delta} \Psi(x,t)}{\partial t^{\delta}}$} \\
&\text{and $\phi_{i}[\Psi(x,t)]=\frac{1}{\omega_{3}}\Psi^2 \Psi^{*}$.}
\end{split}
\end{equation}
Let us now write some of the Adomian polynomials, $\phi_{i}'s$ using the formula in (22) as follows:
\begin{equation*}
\begin{split}
&\text{$\phi_{0}=\frac{1}{\omega_{3}}\Psi_{0}^{2}\Psi_{0}^{*}$,} \\
&\text{$\phi_{1}=\frac{2}{\omega_{3}}\Psi_{0}\Psi_{1}\Psi_{0}^{*}+\frac{1}{\omega_{3}}\Psi_{0}^{2}\Psi_{1}^{*}$,} \\
&\text{$\phi_{2}=\frac{2}{\omega_{3}}\Psi_{0}\Psi_{2}\Psi_{0}^{*}+\frac{1}{\omega_{3}}\Psi_{1}^{2}\Psi_{0}^{*}+
\frac{2}{\omega_{3}}\Psi_{0}\Psi_{1}\Psi_{1}^{*}+\frac{1}{\omega_{3}}\Psi_{0}^{2}\Psi_{2}^{*}$.}
\end{split}
\end{equation*}
By applying the inverse double Laplace transform to the left-hand and right-hand sides of Equation (23), we obtain the following general solution to Equation (16) recursively:
\begin{equation}
\begin{split}
&\Psi_{0}(x,t)=b_{0}(t)+xb_{1}(t), \\
&\Psi_{1}(x,t)= \\
&-(\ell^{x})^{-1}(\ell^{t})^{-1}\displaystyle \left[\frac{1}{s_{1}^{2 \gamma}}\displaystyle \left[\ell^{x}\ell^{t}\displaystyle \left[\frac{\omega_{2}}{\omega_{3}} \frac{\partial^{2 \delta} \Psi_{0}(x,t)}{\partial t^{2 \delta}}+\frac{i}{\omega_{3}} \frac{\partial^{\gamma}\Psi_{0}(x,t)}{\partial x^{\gamma}}+i\frac{\omega_{1}}{\omega_{3}}\frac{\partial^{\delta} \Psi_{0}(x,t)}{\partial t^{\delta}}\right]+
\ell^{x}\ell^{t}[\phi_{0}(\Psi(x,t))]\right]\right], \\
&\text{.} \\
&\text{.} \\
&\text{.} \\
&\Psi_{i+1}(x,t)=-(\ell^{x})^{-1}(\ell^{t})^{-1}\displaystyle \left[\frac{1}{s_{1}^{2 \gamma}}\ell^{x}\ell^{t}\displaystyle \left[R[\Psi_{i}(x,t)]\right]+
\ell^{x}\ell^{t}[\phi_{i}(\Psi(x,t))]\right], \text{for $i\geq0$}.
\end{split}
\end{equation}
Similarly, we can apply the CmDLTr (definition 8) to both sides of Equation (17), we have:
\begin{equation}
\begin{split}
&\displaystyle \ell^{x}\ell^{t}\displaystyle \left[D_{x}^{2\gamma}\Psi\displaystyle \left(\frac{x^{\gamma}}{\gamma},\frac{t^{\delta}}{\delta}\right)\right]=
\tilde{\Psi}(s_{1},s_{2})=\frac{s_{1}^{2 \gamma -1}}{s_{1}^{2 \gamma}}\ell^{x}\ell^{t}[m_{0}(t)]+\frac{s_{1}^{s \gamma -2}}{s_{1}^{2 \gamma}}\ell^{x}\ell^{t}\displaystyle \left[m_{1}\displaystyle \left(\frac{t^{\delta}}{\delta}\right)\right] \\
&-\frac{1}{s_{1}^{2 \gamma}}\displaystyle \left[\ell^{x}\ell^{t} \displaystyle \left[\frac{\omega_{2}}{\omega_{3}} D_{t}^{2\delta}\Psi\displaystyle \left(\frac{x^{\gamma}}{\gamma},\frac{t^{\delta}}{\delta}\right)+\frac{i}{\omega_{3}}
D_{x}^{\gamma}\Psi\displaystyle \left(\frac{x^{\gamma}}{\gamma},\frac{t^{\delta}}{\delta}\right)
+i\frac{\omega_{1}}{\omega_{3}}D_{t}^{\delta}\Psi\displaystyle \left(\frac{x^{\gamma}}{\gamma},\frac{t^{\delta}}{\delta}\right)\right]+
\ell^{x}\ell^{t}\displaystyle \left[\frac{1}{\omega_{3}}\displaystyle \left|\Psi \right|^2 \Psi\right]\right]. \\
&\text{After simplifications, we have:} \\
&\displaystyle \ell^{x}\ell^{t}\displaystyle \left[D_{x}^{2\gamma}\Psi\displaystyle \left(\frac{x^{\gamma}}{\gamma},\frac{t^{\delta}}{\delta}\right)\right]=
\tilde{\Psi}(s_{1},s_{2})=\frac{1}{s_{1}}\ell^{x}\ell^{t}[m_{0}(t)]+\frac{1}{s_{1}^{2}}\ell^{x}\ell^{t}\displaystyle \left[m_{1}\displaystyle \left(\frac{t^{\delta}}{\delta}\right)\right] \\
&-\frac{1}{s_{1}^{2 \gamma}}\displaystyle \left[\ell^{x}\ell^{t} \displaystyle \left[\frac{\omega_{2}}{\omega_{3}} D_{t}^{2\delta}\Psi\displaystyle \left(\frac{x^{\gamma}}{\gamma},\frac{t^{\delta}}{\delta}\right)+\frac{i}{\omega_{3}} D_{x}^{\gamma}\Psi\displaystyle \left(\frac{x^{\gamma}}{\gamma},\frac{t^{\delta}}{\delta}\right)+i\frac{\omega_{1}}{\omega_{3}}
D_{t}^{\delta}\Psi\displaystyle \left(\frac{x^{\gamma}}{\gamma},\frac{t^{\delta}}{\delta}\right)\right]+\ell^{x}\ell^{t}\displaystyle \left[\frac{1}{\omega_{3}}\displaystyle \left[\Psi^2 \Psi^{*}\right]\right]\right].\\
&\text{where $\displaystyle \left|\Psi \right|^2 \Psi=\Psi^2 \Psi^{*}$ such that $\Psi^{*}$ is the conjugate of $\Psi$.}
\end{split}
\end{equation}
Let us now apply the standard NPDEs operator form and (21) to Equation (25), we have:
\begin{equation}
\begin{split}
&\displaystyle \ell^{x}\ell^{t}\displaystyle \left[\sum_{i=0}^{\infty}\Psi_{i}(x,t)\right]=
\frac{1}{s_{1}}\ell^{x}\ell^{t}[m_{0}(t)]+\frac{1}{s_{1}^{2}}\ell^{x}\ell^{t}\displaystyle \left[m_{1}\displaystyle \left(\frac{t^{\delta}}{\delta}\right)\right] \\
&-\frac{1}{s_{1}^{2 \gamma}}\displaystyle \left[\ell^{x}\ell^{t} \displaystyle \left[ R \displaystyle \left[\sum_{i=0}^{\infty}\Psi_{i}(x,t)\right]\right]+\ell^{x}\ell^{t}\displaystyle \left[\sum_{i=0}^{\infty}\phi_{i}    \right]\right]. \\
&\text{where $R[\Psi(x,t)]=\displaystyle \left[\frac{\omega_{2}}{\omega_{3}} D_{t}^{2\delta}\Psi\displaystyle \left(\frac{x^{\gamma}}{\gamma},\frac{t^{\delta}}{\delta}\right)+\frac{i}{\omega_{3}} D_{x}^{\gamma}\Psi\displaystyle \left(\frac{x^{\gamma}}{\gamma},\frac{t^{\delta}}{\delta}\right)+i\frac{\omega_{1}}{\omega_{3}}
D_{t}^{\delta}\Psi\displaystyle \left(\frac{x^{\gamma}}{\gamma},\frac{t^{\delta}}{\delta}\right)\right]$} \\
&\text{and $\phi_{i}[\Psi(x,t)]=\frac{1}{\omega_{3}}\Psi^2 \Psi^{*}$.}
\end{split}
\end{equation}
We apply the inverse double Laplace transform to both sides of Equation (26) to obtain the general solution to Equation (17) recursively as follows:
\begin{equation}
\begin{split}
&\Psi_{0}(x,t)=m_{0}\displaystyle \left(\frac{t^{\delta}}{\delta}\right)+\frac{x^{\gamma}}{\gamma}m_{1}\displaystyle \left(\frac{t^{\delta}}{\delta}\right), \\
&\Psi_{1}(x,t)=-(\ell^{x})^{-1}(\ell^{t})^{-1} \times \\
&\displaystyle \left[\frac{1}{s_{1}^{2 \gamma}}\displaystyle \left[\ell^{x}\ell^{t}\displaystyle \left[\frac{\omega_{2}}{\omega_{3}} D_{t}^{2\delta}\Psi_{0}\displaystyle \left(\frac{x^{\gamma}}{\gamma},\frac{t^{\delta}}{\delta}\right)+\frac{i}{\omega_{3}} D_{x}^{\gamma}\Psi_{0}\displaystyle \left(\frac{x^{\gamma}}{\gamma},\frac{t^{\delta}}{\delta}\right)+i\frac{\omega_{1}}{\omega_{3}}
D_{t}^{\delta}\Psi_{0}\displaystyle \left(\frac{x^{\gamma}}{\gamma},\frac{t^{\delta}}{\delta}\right)\right]+
\ell^{x}\ell^{t}[\phi_{0}(\Psi(x,t))]\right]\right], \\
\vdots \\
&\Psi_{i+1}(x,t)=-(\ell^{x})^{-1}(\ell^{t})^{-1}\displaystyle \left[\frac{1}{s_{1}^{2 \gamma}}\ell^{x}\ell^{t}\displaystyle \left[R[\Psi_{i}(x,t)]\right]+
\ell^{x}\ell^{t}[\phi_{i}(\Psi(x,t))]\right], \text{for $i\geq0$}.
\end{split}
\end{equation}
\subsection*{Numerical Experiment 1:} By applying definitions and properties of Caputo fractional derivative and double Laplace transform, the following numerical experiment will solve Equation (16) analytically:
Let $\omega_{1}=\omega_{2}=\omega_{3}=1$, and $b_{0}(t)=e^{it};b_{1}(t)=0;a_{0}(x)=a_{1}(x)=0$ in (16), we have:
\begin{equation}
\begin{split}
&\frac{\partial^{2\gamma} \Psi(x,t)}{\partial x^{2\gamma}}=-\frac{\partial^{2\delta} \Psi(x,t)}{\partial t^{2\delta}}-i\frac{\partial^{\gamma} \Psi(x,t)}{\partial x^{\gamma}}-i\frac{\partial^{\delta} \Psi(x,t)}{\partial t^{\delta}}-\displaystyle \left|\Psi \right|^2 \Psi. \\
&\text{subject to the following initial and boundary conditions:} \\
&\text{$\Psi(x,0)=0$ and $\frac{\partial\Psi(x,0)}{\partial t}=0$}. \\
&\text{$\Psi(0,t)=e^{it}$ and $\frac{\partial\Psi(0,t)}{\partial x}=0$}. \\
&\text{where $i=\sqrt{-1}$, $0< \gamma,\delta \leq 1$, $t,x>0$; $x,t \in\Re^{+}$}.
\end{split}
\end{equation}
To solve Equation (28), we use our result in (24) as follows:
\begin{equation*}
\begin{split}
&\Psi_{0}(x,t)=e^{it}, \\
&\Psi_{1}(x,t)= \\
&-(\ell^{x})^{-1}(\ell^{t})^{-1}\displaystyle \left[\frac{1}{s_{1}^{2 \gamma}}\displaystyle \left[\ell^{x}\ell^{t}\displaystyle \left[\frac{\partial^{2 \delta} \Psi_{0}(x,t)}{\partial t^{2 \delta}}+i\frac{\partial^{\gamma}\Psi_{0}(x,t)}{\partial x^{\gamma}}+i\frac{\partial^{\delta} \Psi_{0}(x,t)}{\partial t^{\delta}}\right]+
\ell^{x}\ell^{t}[\phi_{0}(\Psi(x,t))]\right]\right] \\
&=-(\ell^{x})^{-1}(\ell^{t})^{-1}\displaystyle \left[\frac{1}{s_{1}^{2 \gamma}}\displaystyle \left[\ell^{x}\ell^{t}\displaystyle \left[\frac{\partial^{2 \delta} \Psi_{0}(x,t)}{\partial t^{2 \delta}}+i\frac{\partial^{\gamma}\Psi_{0}(x,t)}{\partial x^{\gamma}}+i\frac{\partial^{\delta} \Psi_{0}(x,t)}{\partial t^{\delta}}\right]+
\ell^{x}\ell^{t}[\Psi_{0}^{2}\Psi_{0}^{*}]\right]\right] \\
&=-\frac{x^{2 \gamma}}{\Gamma(2 \gamma+1)}\displaystyle \left[\displaystyle \left[it^{1-2\delta}E_{1,2-2\delta}(it)-t^{1-\delta}E_{1,2-\delta}(it)\right]+e^{it}\right] \\
&=-\frac{x^{2 \gamma}}{\Gamma(2 \gamma+1)}\displaystyle \left[\displaystyle \left[iE(t,1-2\delta,i)-E(t,1-\delta,i)\right]+e^{it}\right],
\end{split}
\end{equation*}
\begin{equation*}
\begin{split}
&\Psi_{2}(x,t)= \\
&-(\ell^{x})^{-1}(\ell^{t})^{-1}\displaystyle \left[\frac{1}{s_{1}^{2 \gamma}}\displaystyle \left[\ell^{x}\ell^{t}\displaystyle \left[\frac{\partial^{2 \delta} \Psi_{1}(x,t)}{\partial t^{2 \delta}}+i\frac{\partial^{\gamma}\Psi_{1}(x,t)}{\partial x^{\gamma}}+i\frac{\partial^{\delta} \Psi_{1}(x,t)}{\partial t^{\delta}}\right]+
\ell^{x}\ell^{t}[\phi_{1}(\Psi(x,t))]\right]\right] \\
&=-(\ell^{x})^{-1}(\ell^{t})^{-1}\displaystyle \left[\frac{1}{s_{1}^{2 \gamma}}\displaystyle \left[\ell^{x}\ell^{t}\displaystyle \left[\frac{\partial^{2 \delta} \Psi_{1}(x,t)}{\partial t^{2 \delta}}+i\frac{\partial^{\gamma}\Psi_{1}(x,t)}{\partial x^{\gamma}}+i\frac{\partial^{\delta} \Psi_{1}(x,t)}{\partial t^{\delta}}\right]+
\ell^{x}\ell^{t}[2\Psi_{0}\Psi_{1}\Psi_{0}^{*}+\Psi_{0}^{2}\Psi_{1}^{*}]\right]\right] \\
&=-\frac{x^{2 \gamma}}{\Gamma(2 \gamma+1)}\displaystyle \left[-\frac{x^{2 \gamma}}{\Gamma(2 \gamma+1)}\frac{\partial^{2\delta}}{\partial t^{2\delta}}\displaystyle \left(\displaystyle \left[iE(t,1-2\delta,i)-E(t,1-\delta,i)\right]+e^{it}\right)\right] \\
&-i\frac{x^{2 \gamma}}{\Gamma(2 \gamma+1)}\displaystyle \left[\displaystyle \left[\left(\displaystyle \left[iE(t,1-2\delta,i)-E(t,1-\delta,i)\right]+e^{it}\right)\right]\frac{\partial^{\gamma}}{\partial x^{\gamma}}\displaystyle \left(-\frac{x^{2 \gamma}}{\Gamma(2 \gamma+1)}\right)\right] \\
&-i\frac{x^{2 \gamma}}{\Gamma(2 \gamma+1)}\displaystyle \left[-\frac{x^{2 \gamma}}{\Gamma(2 \gamma+1)}\frac{\partial^{\delta}}{\partial t^{\delta}}\displaystyle \left(\displaystyle \left[iE(t,1-2\delta,i)-E(t,1-\delta,i)\right]+e^{it}\right)\right] \\
&-\frac{x^{2 \gamma}}{\Gamma(2 \gamma+1)}\displaystyle \left[(2e^{it})\displaystyle \left(-\frac{x^{2 \gamma}}{\Gamma(2 \gamma+1)}\right)\displaystyle \left(\displaystyle \left[iE(t,1-2\delta,i)-E(t,1-\delta,i)\right]+e^{it}\right)(e^{-it})\right] \\
&-\frac{x^{2 \gamma}}{\Gamma(2 \gamma+1)}\displaystyle \left[(e^{2it})\displaystyle \left(-\frac{x^{2 \gamma}}{\Gamma(2 \gamma+1)}\right)\displaystyle \left(\displaystyle \left[-iE(t,1-2\delta,i)-E(t,1-\delta,i)\right]+e^{-it}\right)(e^{-it})\right]
\end{split}
\end{equation*}
\begin{equation*}
\begin{split}
&=\frac{x^{4 \gamma}}{\Gamma(4 \gamma+1)}\displaystyle \left[iE(t,1-4\delta,i)-E(t,1-3\delta,i)+iE(t,1-2\delta,i)      \right] \\
&+\frac{x^{3 \gamma}}{\Gamma(3 \gamma+1)}\displaystyle \left[E(t,1-2\delta,i)-iE(t,1-\delta,i)+ie^{it}\right] \\
&+\frac{x^{4 \gamma}}{\Gamma(4 \gamma+1)}\displaystyle \left[-E(t,1-3\delta,i)-iE(t,1-2\delta,i)-E(t,1-\delta,i)\right] \\
&+\frac{x^{4 \gamma}}{\Gamma(4 \gamma+1)}\displaystyle \left[-2iE(t,1-2\delta,i)-E(t,1-\delta,i)+e^{it}\right] \\
&+\frac{x^{4 \gamma}}{\Gamma(4 \gamma+1)}\displaystyle \left[-ie^{2it}E(t,1-2\delta,i)-e^{2it}E(t,1-\delta,i)+e^{it}\right] \\
&=\frac{x^{3 \gamma}}{\Gamma(3 \gamma+1)}\displaystyle \left[E(t,1-2\delta,i)-iE(t,1-\delta,i)+ie^{it}\right] \\
&+\frac{x^{4 \gamma}}{\Gamma(4 \gamma+1)}\displaystyle \{iE(t,1-4\delta,i)-E(t,1-3\delta,i)+iE(t,1-2\delta,i)
-E(t,1-3\delta,i) \\
&-iE(t,1-2\delta,i)-E(t,1-\delta,i)-2iE(t,1-2\delta,i)-E(t,1-\delta,i)+e^{it} \\
&-ie^{2it}E(t,1-2\delta,i)-e^{2it}E(t,1-\delta,i)+e^{it} \}.
\end{split}
\end{equation*}
\begin{equation*}
\begin{split}
\vdots \\
&\text{and so on.}
\end{split}
\end{equation*}
By using all above obtained results, the general approximate-analytical solution to Equation (28) can be written as follows:
\begin{equation}
\begin{split}
&\Psi(x,t)=e^{it}-\frac{x^{2 \gamma}}{\Gamma(2 \gamma+1)}\displaystyle \left[\displaystyle \left[iE(t,1-2\delta,i)-E(t,1-\delta,i)\right]+e^{it}\right] \\
&+\frac{x^{3 \gamma}}{\Gamma(3 \gamma+1)}\displaystyle \left[E(t,1-2\delta,i)-iE(t,1-\delta,i)+ie^{it}\right]+\frac{x^{4 \gamma}}{\Gamma(4 \gamma+1)}\displaystyle \{iE(t,1-4\delta,i) \\
&-E(t,1-3\delta,i)+iE(t,1-2\delta,i)-E(t,1-3\delta,i)-iE(t,1-2\delta,i)-E(t,1-\delta,i) \\
&-2iE(t,1-2\delta,i)-E(t,1-\delta,i)+e^{it}-ie^{2it}E(t,1-2\delta,i)-e^{2it}E(t,1-\delta,i)+e^{it} \}. \\
&+\ldots
\end{split}
\end{equation}
Hence, the approximate-analytical solution for the MNLSE in (1) in the sense of Caputo fractional derivative has been easily obtained via the double Laplace transform coupled with the Adomian decomposition method.
\subsection*{Numerical Experiment 2:} By applying definitions and properties of conformable derivative in \cite{SiLa, AJ} and double Laplace transform, the following numerical experiment will solve Equation (17) analytically:
Let $\omega_{1}=\omega_{2}=\omega_{3}=1$, and $m_{0}(t)=e^{i\frac{t^{\delta}}{\delta}};m_{1}(t)=0;n_{0}(x)=n_{1}(x)=0$ in (17), we have:
\begin{equation}
\begin{split}
&D_{x}^{2\gamma}\Psi\displaystyle \left(\frac{x^{\gamma}}{\gamma},\frac{t^{\delta}}{\delta}\right)=- D_{t}^{2\delta}\Psi\displaystyle \left(\frac{x^{\gamma}}{\gamma},\frac{t^{\delta}}{\delta}\right)-i D_{x}^{\gamma}\Psi\displaystyle \left(\frac{x^{\gamma}}{\gamma},\frac{t^{\delta}}{\delta}\right)-i D_{t}^{\delta}\Psi\displaystyle \left(\frac{x^{\gamma}}{\gamma},\frac{t^{\delta}}{\delta}\right) \\
&-\displaystyle \left|\Psi \right|^2 \Psi. \\
&\text{subject to the following initial and boundary conditions:} \\
&\text{$\Psi\displaystyle \left(\frac{x^{\gamma}}{\gamma},0\right)=0$ and $D_{t}\Psi\displaystyle \left(\frac{x^{\gamma}}{\gamma},0\right)=0$.} \\
&\text{$\Psi\displaystyle \left(0,\frac{t^{\delta}}{\delta}\right)=e^{i\frac{t^{\delta}}{\delta}}$ and $D_{x}\Psi\displaystyle \left(0,\frac{t^{\delta}}{\delta}\right)=0$.} \\
&\text{where $i=\sqrt{-1}$, $0< \gamma,\delta \leq 1$, $t,x>0$; $x,t \in\Re^{+}$.} \\
\end{split}
\end{equation}
To solve Equation (30), we use our result in (27) as follows:
\begin{equation*}
\begin{split}
&\Psi_{0}(x,t)=e^{i\frac{t^{\delta}}{\delta}}, \\
&\Psi_{1}(x,t)=-(\ell^{x})^{-1}(\ell^{t})^{-1} \times \\
&\displaystyle \left[\frac{1}{s_{1}^{2 \gamma}}\displaystyle \left[\ell^{x}\ell^{t}\displaystyle \left[ D_{t}^{2\delta}\Psi_{0}\displaystyle \left(\frac{x^{\gamma}}{\gamma},\frac{t^{\delta}}{\delta}\right)+i D_{x}^{\gamma}\Psi_{0}\displaystyle \left(\frac{x^{\gamma}}{\gamma},\frac{t^{\delta}}{\delta}\right)+i
D_{t}^{\delta}\Psi_{0}\displaystyle \left(\frac{x^{\gamma}}{\gamma},\frac{t^{\delta}}{\delta}\right)\right]+
\ell^{x}\ell^{t}[\phi_{0}(\Psi(x,t))]\right]\right] \\
&=-(\ell^{x})^{-1}(\ell^{t})^{-1}
\displaystyle \left[\frac{1}{s_{1}^{2 \gamma}}\displaystyle \left[\ell^{x}\ell^{t}\displaystyle \left[ D_{t}^{2\delta}\Psi_{0}\displaystyle \left(\frac{x^{\gamma}}{\gamma},\frac{t^{\delta}}{\delta}\right)+i D_{x}^{\gamma}\Psi_{0}\displaystyle \left(\frac{x^{\gamma}}{\gamma},\frac{t^{\delta}}{\delta}\right)+i
D_{t}^{\delta}\Psi_{0}\displaystyle \left(\frac{x^{\gamma}}{\gamma},\frac{t^{\delta}}{\delta}\right)\right]+
\ell^{x}\ell^{t}[\Psi_{0}^{2}\Psi_{0}^{*}]\right]\right] \\
&=-\frac{x^{2 \gamma-1}}{\gamma^{2\gamma -1}\Gamma(2 \gamma)}\displaystyle \left[\displaystyle \left[-e^{i\frac{t^{\delta}}{\delta}}-e^{i\frac{t^{\delta}}{\delta}}\right]+e^{i\frac{t^{\delta}}{\delta}}\right] \\
&=-\frac{x^{2 \gamma-1}}{\gamma^{2\gamma -1}\Gamma(2 \gamma)}\displaystyle \left[-e^{i\frac{t^{\delta}}{\delta}}\right]
\end{split}
\end{equation*}
\begin{equation*}
\begin{split}
&\Psi_{2}(x,t)=-(\ell^{x})^{-1}(\ell^{t})^{-1} \times \\
&\displaystyle \left[\frac{1}{s_{1}^{2 \gamma}}\displaystyle \left[\ell^{x}\ell^{t}\displaystyle \left[ D_{t}^{2\delta}\Psi_{1}\displaystyle \left(\frac{x^{\gamma}}{\gamma},\frac{t^{\delta}}{\delta}\right)+i D_{x}^{\gamma}\Psi_{1}\displaystyle \left(\frac{x^{\gamma}}{\gamma},\frac{t^{\delta}}{\delta}\right)+i
D_{t}^{\delta}\Psi_{1}\displaystyle \left(\frac{x^{\gamma}}{\gamma},\frac{t^{\delta}}{\delta}\right)\right]+
\ell^{x}\ell^{t}[\phi_{1}(\Psi(x,t))]\right]\right] \\
&=-(\ell^{x})^{-1}(\ell^{t})^{-1} \times \\
&\displaystyle \left[\frac{1}{s_{1}^{2 \gamma}}\displaystyle \left[\ell^{x}\ell^{t}\displaystyle \left[ D_{t}^{2\delta}\Psi_{1}\displaystyle \left(\frac{x^{\gamma}}{\gamma},\frac{t^{\delta}}{\delta}\right)+i D_{x}^{\gamma}\Psi_{1}\displaystyle \left(\frac{x^{\gamma}}{\gamma},\frac{t^{\delta}}{\delta}\right)+i
D_{t}^{\delta}\Psi_{1}\displaystyle \left(\frac{x^{\gamma}}{\gamma},\frac{t^{\delta}}{\delta}\right)\right]+
\ell^{x}\ell^{t}[2\Psi_{0}\Psi_{1}\Psi_{0}^{*}+\Psi_{0}^{2}\Psi_{1}^{*}]\right]\right] \\
&=-\frac{x^{2 \gamma-1}}{\gamma^{2\gamma -1}\Gamma(2 \gamma)}\displaystyle \{\displaystyle \left(-\frac{x^{2 \gamma-1}}{\gamma^{2\gamma -1}\Gamma(2 \gamma)}\right)e^{i\frac{t^{\delta}}{\delta}}+(2\gamma-1)x^{\gamma-1}\displaystyle \left(i\frac{e^{i\frac{t^{\delta}}{\delta}}}{\gamma^{2\gamma -1}\Gamma(2 \gamma)}\right)-\displaystyle \left(\frac{x^{2 \gamma-1}}{\gamma^{2\gamma -1}\Gamma(2 \gamma)}\right)e^{i\frac{t^{\delta}}{\delta}} \\
&-\displaystyle \left(\frac{2x^{2 \gamma-1}}{\gamma^{2\gamma -1}\Gamma(2 \gamma)}\right)e^{-i\frac{t^{\delta}}{\delta}}+
\displaystyle \left(\frac{x^{2 \gamma-1}}{\gamma^{2\gamma -1}\Gamma(2 \gamma)}\right)e^{i\frac{t^{\delta}}{\delta}}\} \\
&=\frac{x^{4 \gamma-2}}{\gamma^{4\gamma -2}\Gamma(4 \gamma)}e^{i\frac{t^{\delta}}{\delta}}-
i\frac{(2\gamma-1)x^{3 \gamma-2}}{\gamma^{4\gamma -2}\Gamma(4 \gamma)}e^{i\frac{t^{\delta}}{\delta}}+
\frac{x^{4 \gamma-2}}{\gamma^{4\gamma -2}\Gamma(4 \gamma)}e^{i\frac{t^{\delta}}{\delta}}+
\frac{2x^{4 \gamma-2}}{\gamma^{4\gamma -2}\Gamma(4 \gamma)}e^{-i\frac{t^{\delta}}{\delta}}-
\frac{x^{4 \gamma-2}}{\gamma^{4\gamma -2}\Gamma(4 \gamma)}e^{i\frac{t^{\delta}}{\delta}} \\
&=-i\frac{(2\gamma-1)x^{3 \gamma-2}}{\gamma^{4\gamma -2}\Gamma(4 \gamma)}e^{i\frac{t^{\delta}}{\delta}}+
\frac{x^{4 \gamma-2}}{\gamma^{4\gamma -2}\Gamma(4 \gamma)}\displaystyle \left[e^{i\frac{t^{\delta}}{\delta}}+e^{i\frac{t^{\delta}}{\delta}}+2e^{-i\frac{t^{\delta}}{\delta}}-e^{i\frac{t^{\delta}}{\delta}}\right] \\
&=-i\frac{(2\gamma-1)x^{3 \gamma-2}}{\gamma^{4\gamma -2}\Gamma(4 \gamma)}e^{i\frac{t^{\delta}}{\delta}}+
\frac{x^{4 \gamma-2}}{\gamma^{4\gamma -2}\Gamma(4 \gamma)}\displaystyle \left[e^{i\frac{t^{\delta}}{\delta}}+2e^{-i\frac{t^{\delta}}{\delta}}\right]. \\
&\vdots \\
&\text{and so on.}
\end{split}
\end{equation*}
By using the above obtained results, the general approximate-analytical solution to Equation (30) can be written as follows:
\begin{equation}
\begin{split}
&\Psi(x,t)=e^{i\frac{t^{\delta}}{\delta}}-\frac{x^{2 \gamma-1}}{\gamma^{2\gamma -1}\Gamma(2 \gamma)}\displaystyle \left[-e^{i\frac{t^{\delta}}{\delta}}\right]-i\frac{(2\gamma-1)x^{3 \gamma-2}}{\gamma^{4\gamma -2}\Gamma(4 \gamma)}e^{i\frac{t^{\delta}}{\delta}}+
\frac{x^{4 \gamma-2}}{\gamma^{4\gamma -2}\Gamma(4 \gamma)}\displaystyle \left[e^{i\frac{t^{\delta}}{\delta}}+2e^{-i\frac{t^{\delta}}{\delta}}\right] \\
&+\ldots
\end{split}
\end{equation}
Hence, the approximate-analytical solution for the MNLSE in (1) in the sense of conformable derivative has also been easily obtained via the double Laplace transform coupled with the Adomian decomposition method.

\section{The graphical comparisons of solutions}
\begin{figure}[htb!]
  \includegraphics[width=\linewidth]{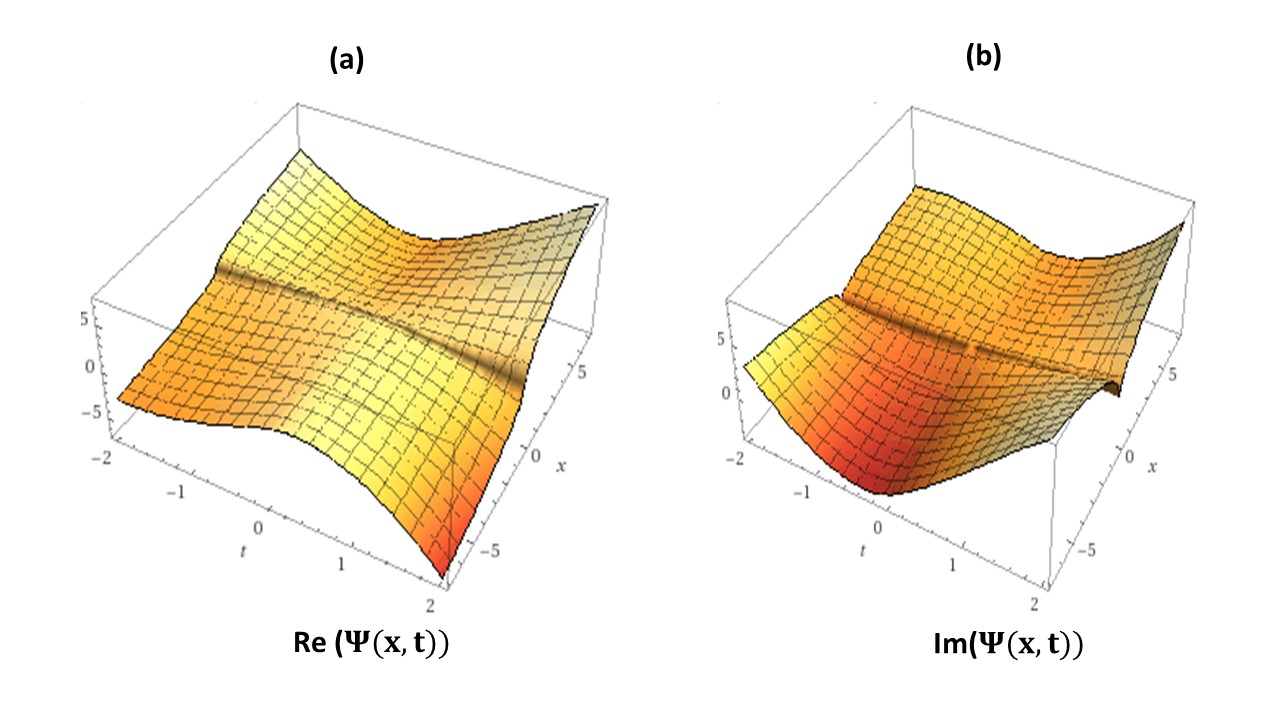}
  \caption{3D Plot of the real part (a) and imaginary part (b) of the Approximate Analytical Solution in (29) for $\gamma=\delta=0.25$}\label{fig:Slide1}
\end{figure}
\begin{figure}[htb!]
  \includegraphics[width=\linewidth]{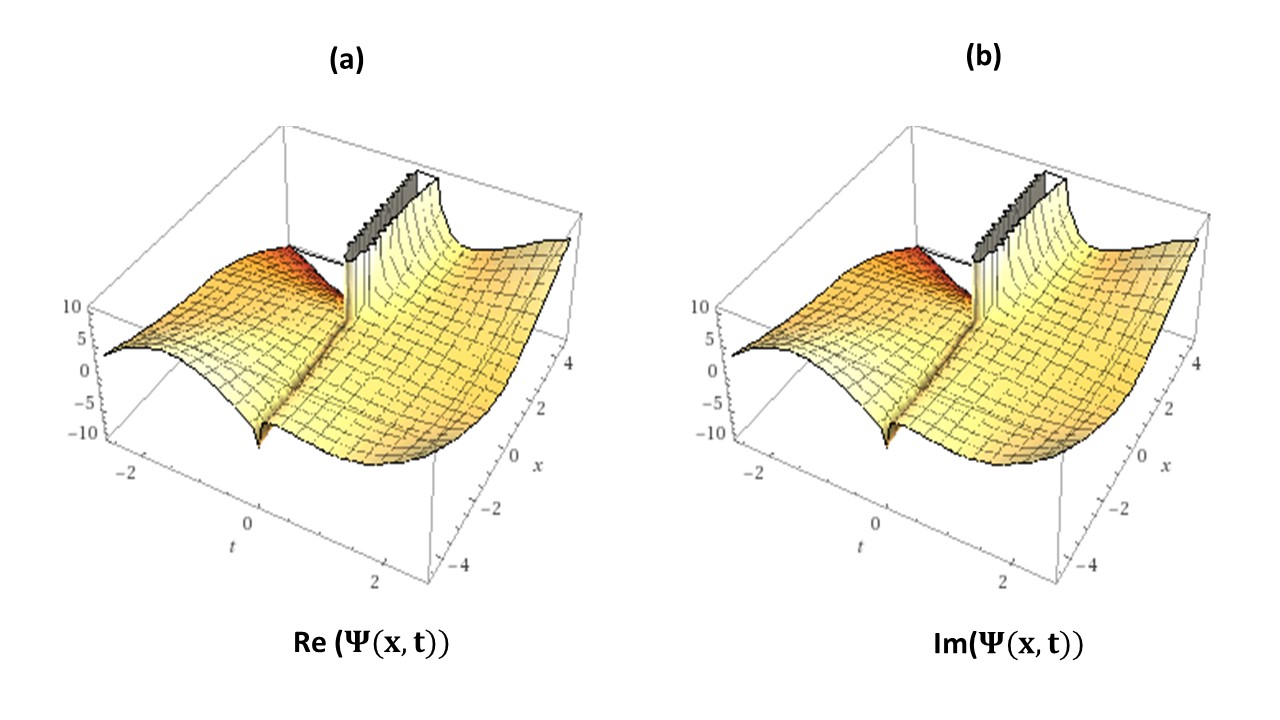}
  \caption{3D Plot of the real part (a) and imaginary part (b) of the Approximate Analytical Solution in (29) for $\gamma=\delta=0.75$}\label{fig:Slide2}
\end{figure}
\begin{figure}[htb!]
  \includegraphics[width=\linewidth]{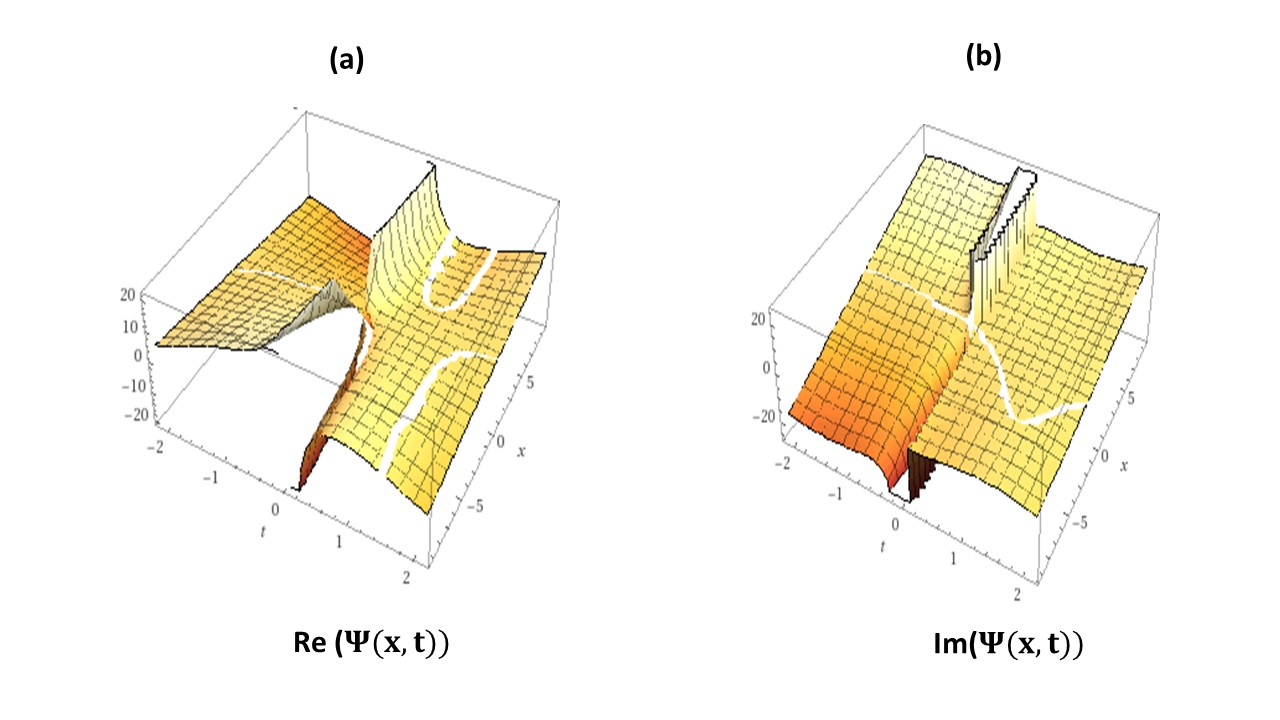}
  \caption{3D Plot of the real part (a) and imaginary part (b) of the Approximate Analytical Solution in (29) for $\gamma=0.50;\delta=0.85$}\label{fig:Slide3}
\end{figure}
\begin{figure}[htb!]
  \includegraphics[width=\linewidth]{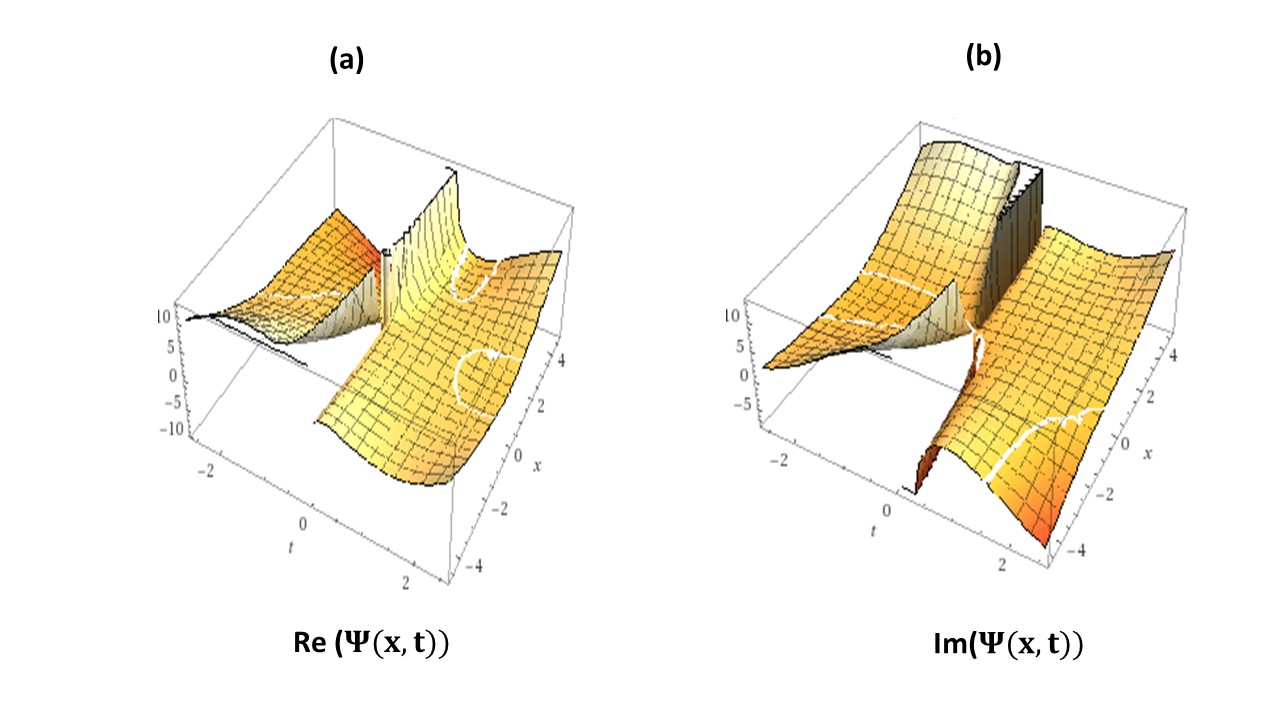}
  \caption{3D Plot of the real part (a) and imaginary part (b) of the Approximate Analytical Solution in (29) for $\gamma=0.75;\delta=0.85$}\label{fig:Slide4}
\end{figure}
\begin{figure}[htb!]
  \includegraphics[width=\linewidth]{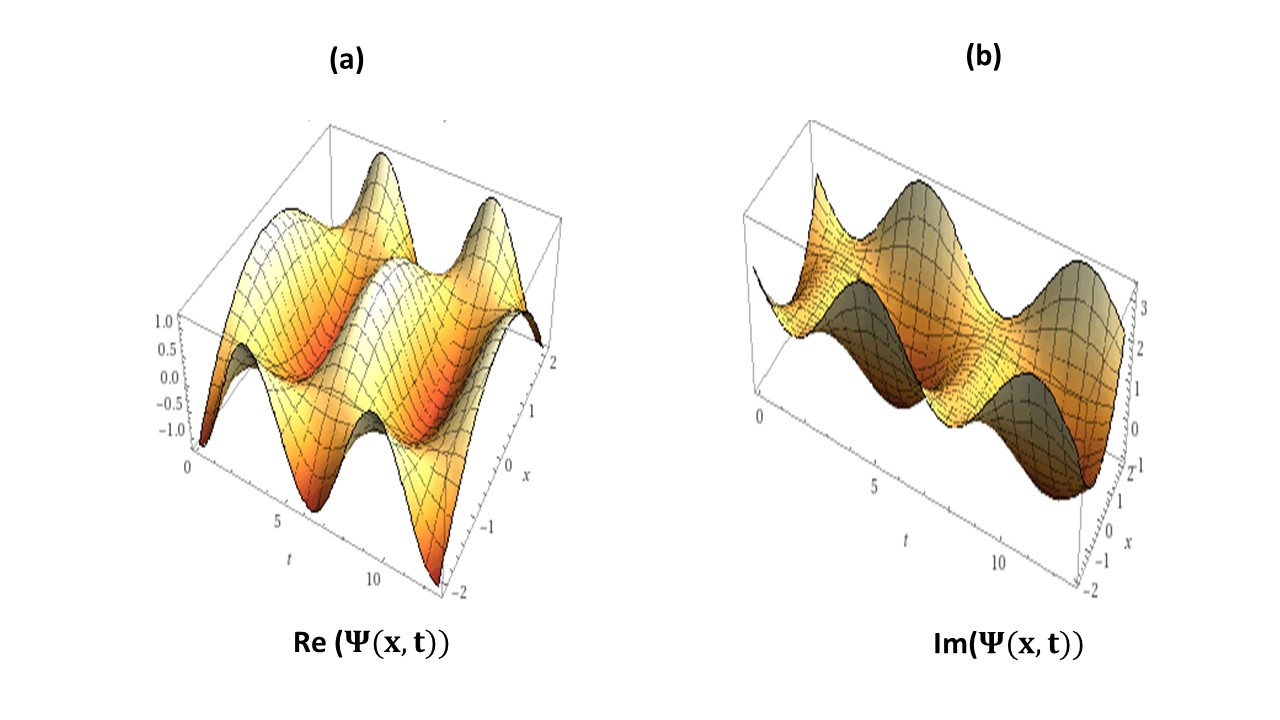}
  \caption{3D Plot of the real part (a) and imaginary part (b) of the Approximate Analytical Solution in (29) for $\gamma=\delta=1$}\label{fig:Slide5}
\end{figure}
\begin{figure}[htb!]
  \includegraphics[width=\linewidth]{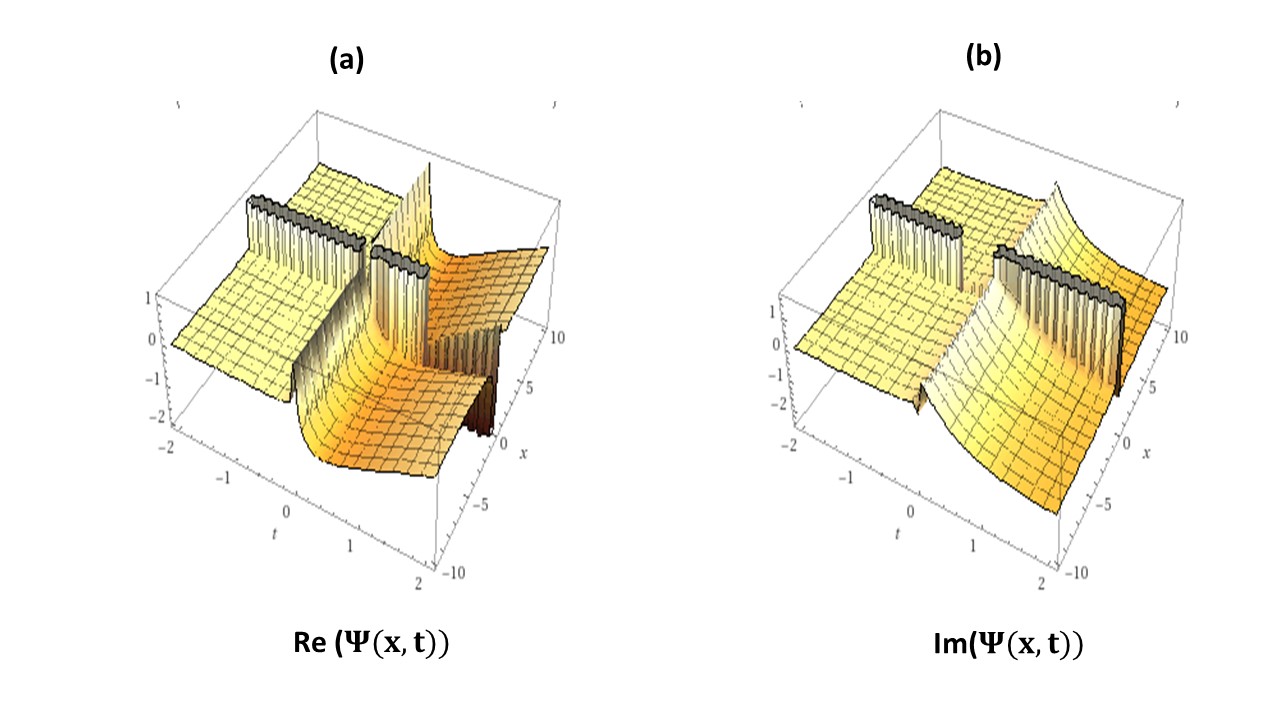}
  \caption{3D Plot of the real part (a) and imaginary part (b) of the Approximate Analytical Solution in (31) for $\gamma=\delta=0.25$}\label{fig:Slide6}
\end{figure}
\begin{figure}[htb!]
  \includegraphics[width=\linewidth]{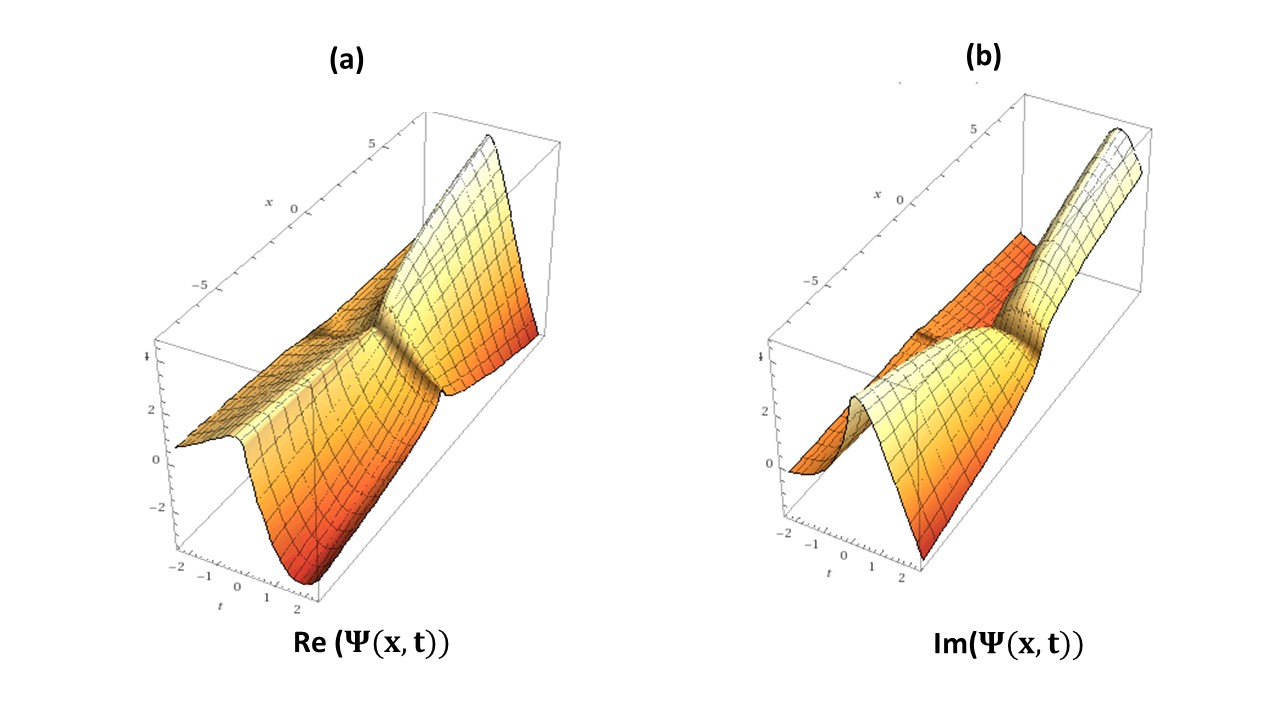}
  \caption{3D Plot of the real part (a) and imaginary part (b) of the Approximate Analytical Solution in (31) for $\gamma=\delta=0.75$}\label{fig:Slide7}
\end{figure}
\begin{figure}[htb!]
  \includegraphics[width=\linewidth]{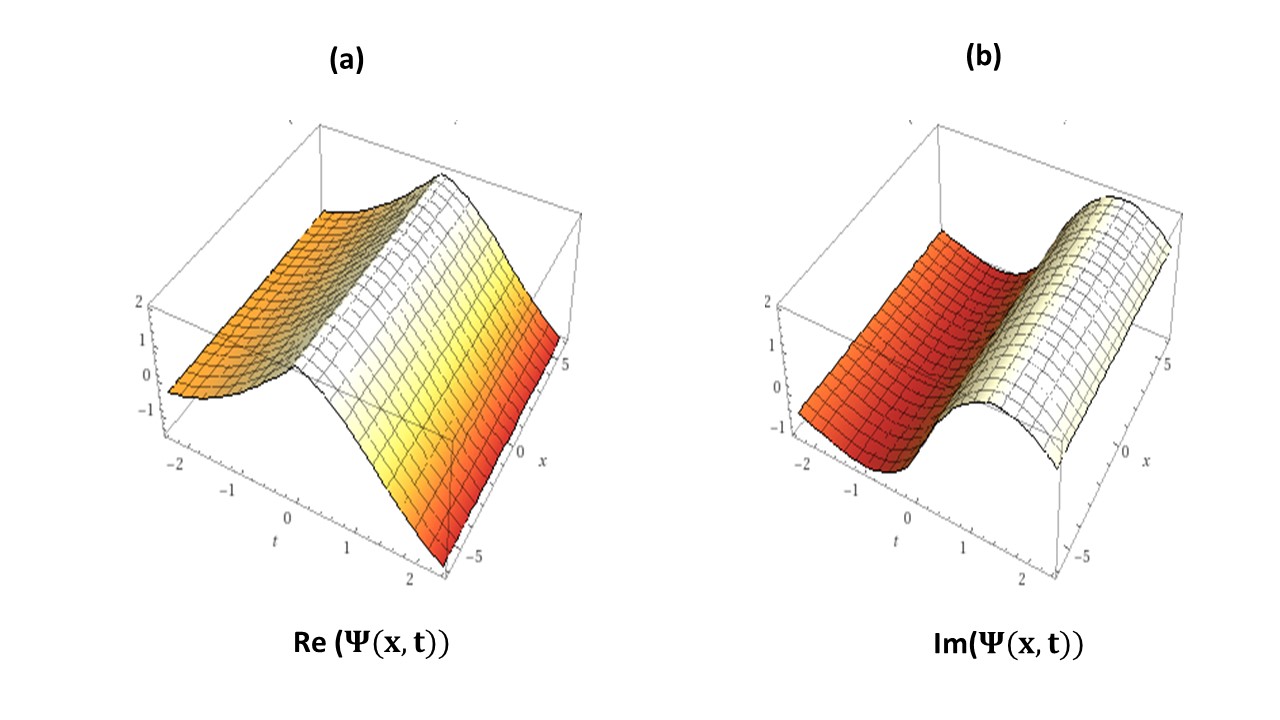}
  \caption{3D Plot of the real part (a) and imaginary part (b) of the Approximate Analytical Solution in (31) for $\gamma=0.50;\delta=0.85$}\label{fig:Slide8}
\end{figure}
\begin{figure}[htb!]
  \includegraphics[width=\linewidth]{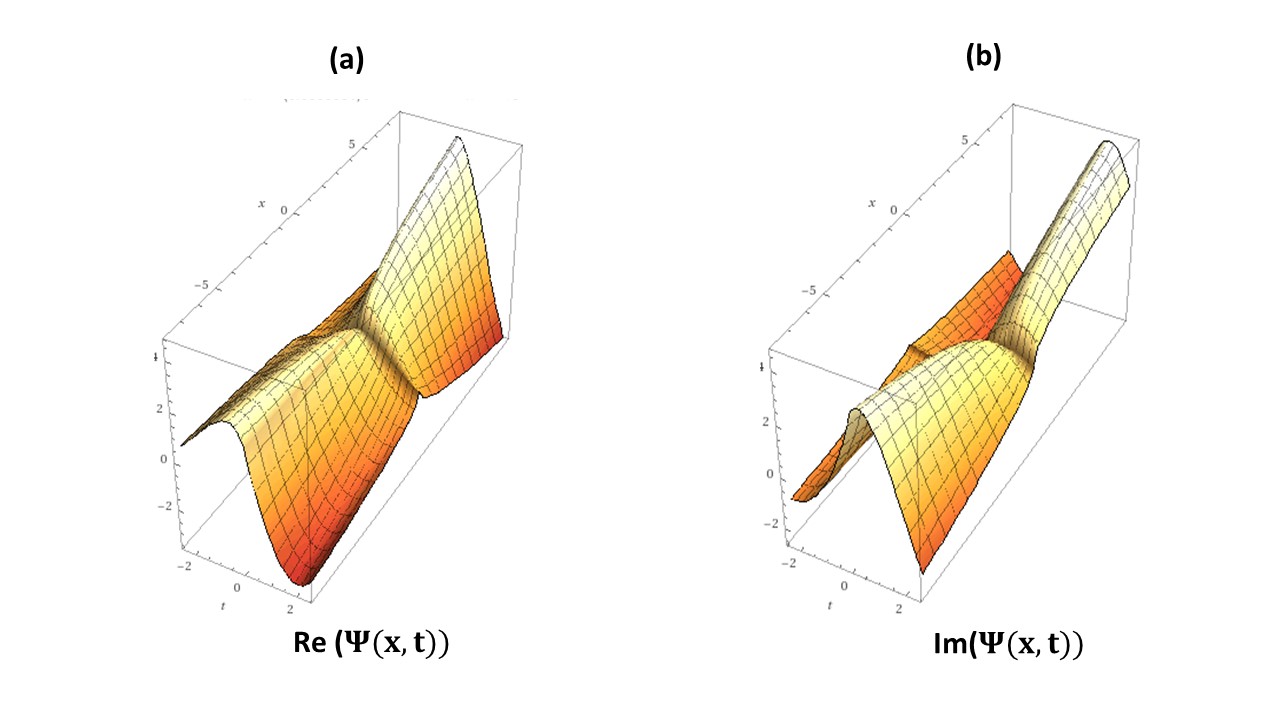}
  \caption{3D Plot of the real part (a) and imaginary part (b) of the Approximate Analytical Solution in (31) for $\gamma=0.75;\delta=0.85$}\label{fig:Slide9}
\end{figure}
\begin{figure}[htb!]
  \includegraphics[width=\linewidth]{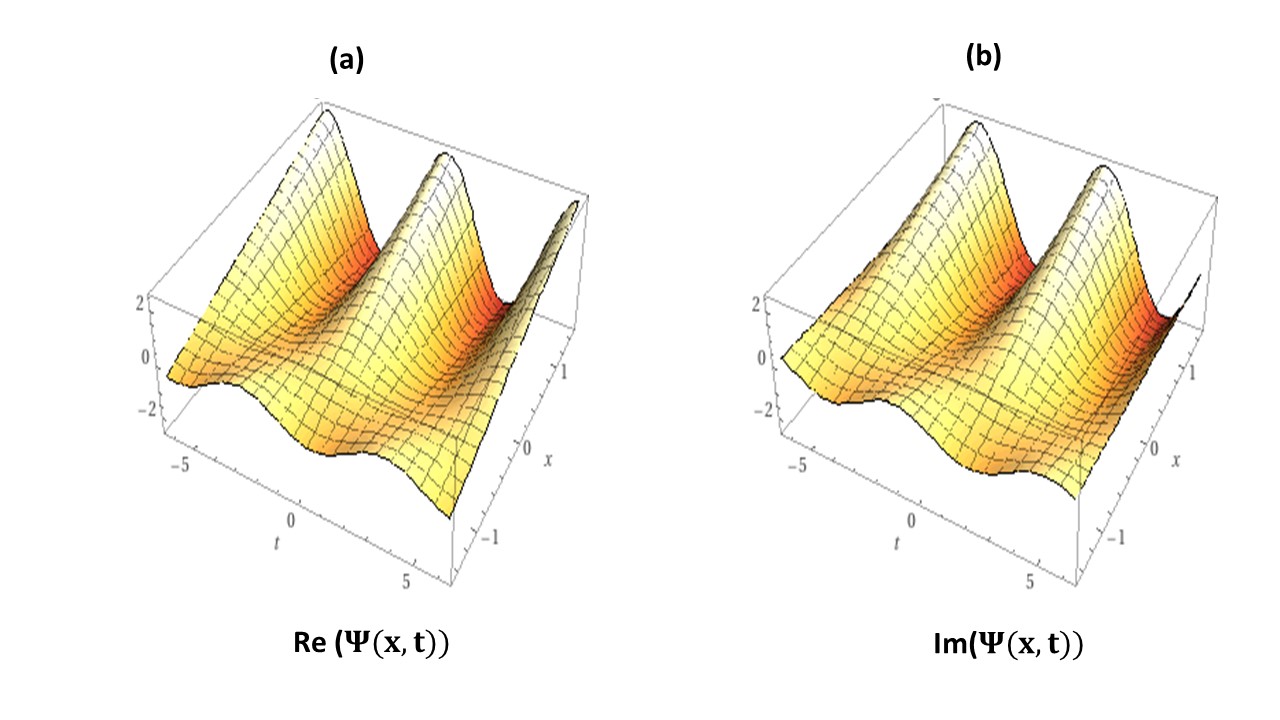}
  \caption{3D Plot of the real part (a) and imaginary part (b) of the Approximate Analytical Solution in (31) for $\gamma=\delta=1$}\label{fig:Slide10}
\end{figure}

\begin{table}[htb!]
\centering
\caption{Comparison of the absolute approximate solutions from CpDLTr, CmDLTr, and the exact solution of the space-time fractional nonlinear Schr\"{o}dinger equation (Example 9) in \cite{CHinD}. Note: $Error_{cp}=|Exact-CpDLTr|$ and $Error_{cm}=|Exact-CmDLTr|$.}
\label{my-label}
\begin{tabular}{|l|l|l|l|l|l|l|l|l|}
\hline
\multirow{2}{*}{$(x,t)$} & \multirow{2}{*}{Exact} & \multirow{2}{*}{$\gamma$} & \multirow{2}{*}{$\delta$} & \multirow{2}{*}{CpDLTr} & \multirow{2}{*}{CmDLTr} & \multirow{2}{*}{$Error_{cp}$} & \multirow{2}{*}{$Error_{cm}$}\\
 &  &      &      &        &      &      &  \\ \hline
\multirow{3}{*}{(0.1,0.1)} & \multirow{3}{*}{0.108060} & 0.25 & 0.25 & 0.034850     & 0.627680 & 0.073309 & 0.519620\\ \cline{3-8}
 &  & 0.75 & 0.75 & 0.052744 & 0.972022  & 0.055315 & 0.863962 \\ \cline{3-8}
 &  & 1    & 1    & 0.060835 & 0.058760  & 0.047224 & 0.049299 \\ \hline
\multirow{3}{*}{(0.3,0.3)} & \multirow{3}{*}{0.324180}  & 0.25 & 0.25 & 0.061912   & 0.983617 & 0.262268 & 0.659436\\ \cline{3-8}
 &  & 0.75 & 0.75 & 0.142055 & 0.857462  & 0.182126 & 0.533281 \\ \cline{3-8}
 &  & 1    & 1    & 0.223340 & 0.199797  & 0.100841 & 0.124384 \\ \hline
\multirow{3}{*}{(0.5,0.5)} & \multirow{3}{*}{0.540300}  & 0.25 & 0.25 & 0.054601   & 0.975461 & 0.485698 & 0.435161\\ \cline{3-8}
 &  & 0.75 & 0.75 & 0.198301 & 0.701852  & 0.341999 & 0.161552 \\ \cline{3-8}
 &  & 1    & 1    & 0.440288 & 0.361382  & 0.100012 & 0.178918 \\ \hline
\multirow{3}{*}{(0.7,0.7)} & \multirow{3}{*}{0.756420}  & 0.25 & 0.25 & 0.021153   & 0.869221 & 0.735269 & 0.112798\\ \cline{3-8}
 &  & 0.75 & 0.75 & 0.211584 & 0.523044  & 0.544839 & 0.233379 \\ \cline{3-8}
 &  & 1    & 1    & 0.711680 & 0.530550  & 0.044743 & 0.225873 \\ \hline
\multirow{3}{*}{(0.9,0.9)} & \multirow{3}{*}{0.972540}  & 0.25 & 0.25 & 0.034267   & 0.728667 & 0.938276 & 0.243877\\ \cline{3-8}
 &  & 0.75 & 0.75 & 0.173940 & 0.332328  & 0.798597 & 0.640216 \\ \cline{3-8}
 &  & 1    & 1    & 1.037520 & 0.694332  & 0.064976 & 0.278212 \\ \hline
\end{tabular}
\end{table}
In this section, the obtained approximate solutions in both (29) and (31) have been graphically compared for various values of $\gamma$ and $\delta$ (see figures 1 to 10 ) where each graph shows both real part and imaginary part of solution. In addition, table 1 provides a comparison of absolute approximate solutions from CpDLTr in (29) and from CmDLTr in (31) at $\gamma=\delta=0.25;\gamma=\delta=0.75;\gamma=\delta=1$ with the exact solution from Example 9 in \cite{CHinD}. According to table 1, at $x=t=0.1;x=t=0.3;x=t=0.5$, the absolute error value from exact and the approximate solution from CpDLTr is less than the absolute error value from exact and the approximate solution from CmDLTr. At $x=t=0.7$, for $\gamma=\delta=0.25;0.75$ the absolute error value from exact and the approximate solution from CmDLTr is less than the absolute error value from exact and the approximate solution from CpDLTr, while at $x=t=0.7$ for $\gamma=\delta=1$, the the approximate solution from CpDLTr converges to the exact solution better than the one from CmDLTr. Similarly, at $x=t=0.9$ for $\gamma=\delta=1$, the approximate solution from CpDLTr converges to the exact solution better than the one from CmDLTr. Therefore, the obtained approximate solution in the sense of Caputo fractional derivative is much better than the obtained solution in the sense of conformable derivatives. On one hand, the Caputo fractional derivative is a nonlocal fractional operator which provides a good interpretation to the physical behavior of systems, while the conformable derivative is a type of local fractional derivative which is basically a generalized form of usual limit-based derivative which lacks some of the important properties to be classified as a fractional derivative. As a result, solving systems of nonlinear partial differential equations in the sense of Caputo fractional derivatives is highly recommended. However, exploring the definition of conformable derivative is also interesting because as the authors believe that any new mathematical definition deserves to be explored and investigated.

\section{Conclusion}

Nonlinear Schr\"{o}dinger equation has been an interesting field of research for many mathematicians and scientists due to the important applications of this equation in physics and engineering. This research study provides a powerful mathematical tool to solve the nonlinear Schr\"{o}dinger equation involving both Caputo fractional derivative and conformable derivative. Therefore, the generalized double Laplace transform method can be efficiently applied in solving nonlinear fractional Schr\"{o}dinger equation and all other nonlinear fractional partial differential equations.

\subsection*{Disclosure statement} The authors declare no conflict of interests.

\subsection*{Acknowledgments} This research did not receive any specific grant from funding agencies in the public, commercial, or not-for-profit sectors.

\end{document}